\documentclass{amsart}

\usepackage{amsbsy,amssymb,amscd,amsfonts,latexsym,amstext,delarray,
amsmath,graphicx} 
\input xypic
\usepackage{color}

\newtheorem{thm}{Theorem}[section]
\newtheorem{prop}[thm]{Proposition}
\newtheorem{cor}[thm]{Corollary}
\newtheorem{lem}[thm]{Lemma}
\newtheorem{conj}{Conjecture}[section]
\newtheorem{defn}[thm]{Definition}
\newtheorem{rem}[thm]{Remark}

\newtheorem{ques}[thm]{Question}

\numberwithin{equation}{section}

\def\bF{{\mathbb F}}
\def\bG{{\mathbb G}}

\def\bK{{\mathbb K}}
\def\bL{{\mathbb L}}

\def\bT{{\mathbb T}}

\def\A{{\mathbb A}}
\def\C{{\mathbb C}}
\def\F{{\mathbb F}}

\renewcommand{\P}{{\mathbb P}}
\def\Q{{\mathbb Q}}
\def\Z{{\mathbb Z}}

\def\cA{{\mathcal A}}

\def\cC{{\mathcal C}}

\def\cE{{\mathcal E}}

\def\cG{{\mathcal G}}
\def\cH{{\mathcal H}}

\def\cL{{\mathcal L}}
\def\cM{{\mathcal M}}
\def\cN{{\mathcal N}}

\def\cT{{\mathcal T}}
\def\cU{{\mathcal U}}
\def\cV{{\mathcal V}}

\def\cX{{\mathcal X}}

\def\m{{\mathfrak{m}}}

\newcommand{\mbf}[1]{\mathbf{#1}}
\newcommand{\csm}{c_{SM}}

\title{Quantum Field Theory over $\bF_1$}
\author{Dori Bejleri and Matilde Marcolli}
\address{Mathematics Department, Mail Code 253-37, Caltech, 1200 E.~California Blvd. Pasadena, CA 91125, USA}
\email{dbejleri@caltech.edu}
\email{matilde@caltech.edu} 

\begin{document}

\begin{abstract}
In this paper we discuss some questions about geometry over the
field with one element, motivated by the properties of algebraic
varieties that arise in perturbative quantum field theory. We follow
the approach to $\F_1$-geometry based on torified schemes. We
first discuss some simple necessary conditions in terms of Euler
characteristic and classes in the Grothendieck ring, then we give
a blowup formula for torified varieties and we show that the
wonderful compactifications of the graph configuration spaces,
that arise in the computation of Feynman integrals in position
space, admit an $\F_1$-structure. By a similar argument we
show that the moduli spaces of curves $\bar\cM_{0,n}$ admit
an $\F_1$-structure. We also discuss conditions on 
hyperplane arrangements, a possible notion of embedded
$\F_1$-structure and its relation to Chern classes,
and questions on Chern classes of varieties with regular torifications.
\end{abstract}

\maketitle

\section{Introduction}

This paper, as the title readily suggests, is inspired by Oliver Schnetz's interesting 
paper ``Quantum field theory over $\F_q$", \cite{Schn}, with the motivation of
investigating under which circumstances one could envision the existence of a
quantum field theory over the field with one element. This is meant not so much
in the sense of developing an actual physical Lagrangian and Feynman rules in the
context of $\F_1$-geometry, but of investigating when certain classes of algebraic
varieties that naturally arise in the context of perturbative quantum field theory,
and which are already defined over $\Z$, may be carrying an additional
$\F_1$-structure. 

\smallskip

The relation between
quantum field theory and geometry of varieties over finite fields originally arises from the
question of ``polynomial countability" for a class of hypersurfaces $X_\Gamma$
defined by the parametric formulation of momentum space 
Feynman integrals in perturbative quantum field theory. 
In addition to these much studied graph hypersurfaces $X_\Gamma$, there are
other algebraic varieties directly connected to the computation of Feynman
integrals. In particular, certain complete intersections $\Lambda_\Gamma$,
recently studied by Esnault and by Bloch, \cite{Bloch}, and some closely related hyperplane
arrangements $\cA_\Gamma$. There are also the graph configuration 
spaces ${\rm Conf}_\Gamma(X)=X^{\# V(\Gamma)}\smallsetminus \cup_e \Delta_e$
and their wonderful compactifications, which arise in the computation of Feynman
integrals in configuration spaces \cite{CeyMa1}, \cite{CeyMa2}. In all of these cases,
the motives of these varieties carry some useful information about the corresponding
Feynman integral computations.

\medskip

In a different direction, the idea of the ``field with one element" and of
the existence of a suitable notion of 
$\mathbb{F}_1$-geometry, first arose from a comment by J.~Tits in a paper from 1956,
where he noted that 
the number of points of $\text{GL}_n(\mathbb{F}_q)$
is a polynomial $N(q) = (q^n - 1)(q^n - q)...(q^n - q^{n-1})$ with a zero of degree 
$r$ at $q = 1$, such that $\lim_{q \to 1} N(q)/(q - 1)^r = n!$, and suggested that one
might interpret the symmetric group $S_n = \text{GL}_n(\mathbb{F}_1)$ 
as an algebraic group in ``characteristic 1." This suggests, more generally, that 
for any scheme $X$ defined over $\mathbb{Z}$ for which the counting function
\begin{equation}\label{NXq}
N_X(q) = \# X_q(\mathbb{F}_q) 
\end{equation}
giving the number of points over $\mathbb{F}_q$
of the reduction of $X$ , and for which there is a limit  
\begin{equation}\label{Nq1}
 \lim_{q \to 1} N_X(q)/(q - 1)^r ,
\end{equation} 
where $r$ is the order of the zero at $q = 1$, one
should interpret this limit as the number of $\mathbb{F}_1$-points, 
$\# X(\mathbb{F}_1)$ of $X$. 
Simple examples include projective space, for which 
$\# \mathbb{P}^{(n-1)}(\mathbb{F}_1) = n$, and Grassmanians, for
which $\# \text{Gr}(k,n)(\mathbb{F}_1) = {n \choose k}$, see \cite{lorscheid1}. 
Several different forms of geometry over $\F_1$ were developed in recent
years (see \cite{LL2} for a general overview).  

\smallskip

There are at present many different versions of geometry over the field
with one element (see \cite{LL2} for a comparative survey). Throughout
this paper we adopt the approach to $\F_1$-geometry developed
by Javier L\'opez Pe\~{n}a and Oliver Lorscheid, based on the
existence of affine torifications.

\medskip

In this paper, we discuss some aspects of the geometry of the
varieties associated to Feynman integrals, from the point of
view of $\F_1$-geometry. This will lead us to formulate a series 
of questions, about these classes of varieties as well as, more 
generally, about $\F_1$-geometry. 

\smallskip

In \S \ref{QFTsec} we recall the original questions about the motivic
properties of the graph hypersurfaces and we show how a simple
modification of the standard derivation of the parametric Feynman
integral of perturbative QFT leads to the occurrence of other 
algebraic varieties, including the mixed Tate complete intersection
$\Lambda_\Gamma$ recently studied by Esnault and by Bloch.
We also recall the definition of the varieties $\overline{\rm Conf}_\Gamma(X)$,
the wonderful compactifications of graph configuration spaces, used in
\cite{CeyMa1}, \cite{CeyMa2} to study Feynman integrals in configuration
spaces.

\smallskip

In \S \ref{EulGrSec} we describe some simple necessary conditions for
the existence of an $\F_1$-structure, based on constraints on the Euler
characteristic and on the class in the Grothendieck ring of varieties and
we give explicit examples of cases where the graph hypersurfaces
satisfy or not these constraints. 

\smallskip

In \S \ref{ConfF1Sec}, we give a blowup formula for torified varieties
and we use it to show that, in the case of $X=\P^D$, the wonderful
compactifications $\overline{\rm Conf}_\Gamma(X)$ admit a regular
affine torification in the sense of \cite{LL}, hence they are $\F_1$-varieties.
By the same technique, we show that the moduli space of curves
$\bar\cM_{0,n}$ is an $\F_1$-variety and that so are the generalizations
$T_{d,n}$, with $T_{1,n}=\bar\cM_{0,n}$ considered in \cite{CheGibKra}.

\smallskip

In \S \ref{HypSec} we identify hyperplane arrangements as an especially interesting
class of varieties about which to investigate when they admit an $\F_1$-structure.
This is especially so, in view of their (conjectural) role 
as ``generators" of mixed Tate motives. We identify an explicit necessary
condition in terms of the coefficients of the characteristic polynomial of the
arrangement and its M\"obius function, as a direct consequence of a recent 
result of Aluffi \cite{Alu3} on the Grothendieck classes of hyperplane arrangements.

\smallskip

In \S \ref{ChernSec} we use another result of Aluffi \cite{Alu6} expressing
the Chern classes of singular varieties in terms of Euler characteristics
of hyperplane sections, to suggest the existence of a notion of embedded
$\F_1$-structures. 

\smallskip

In \S \ref{TorSec} we formulate the question of the existence of
an analog of the Ehlers formula for Chern classes of toric varieties
in the case of varieties that admit a regular affine torification. While
we show that the argument given in \cite{Alu2} for the Ehlers formula
does not directly extend to the regular torified case, we suggest that
there may be a reformulation in terms of regular torifications of the
results of \cite{AluMi} for Chern classes Schubert varieties.

\medskip
\section{Quantum field theory, algebraic varieties and motives}\label{QFTsec}

We recall here some well known facts about the occurrence of algebraic
varieties, motives and periods in perturbative quantum field theory, which
motivate some of the questions that we analyze in the rest of the paper.

\smallskip

In quantum field theory, Feynman graphs parameterize the perturbative expansion of the Green's functions. To each Feynman graph $\Gamma$, there is an associated (typically divergent) integral $U(\Gamma)$ that gives a term in this expansion.  Formally (ignoring divergences)
each of these integrals can be written in the Feynman parametric form as an integral of an
algebraic differential form on the complement of a (singular) projective hypersurface $X_\Gamma
\subset \P^{n-1}$, with $n$ the number of edges of the graph, defined by the vanishing
of a graph polynomial, the Kirchhoff polynomial $\Psi_\Gamma$, see \cite{BEK}, \cite{Mar}. 
Actually accounting for divergences introduces serious complications, which require
performing blowups that separate the locus of integration from the hypersurface and
accounting for the ambiguities deriving from monodromies, see \cite{BlKr}. In several
physically significant cases, it has been shown \cite{Bro} that the resulting period obtained
from the evaluation of the Feynman amplitude $U(\Gamma)$ is a period of a mixed Tate
motive over $\Z$, that is, a $\Q [1/(2\pi i)]$-linear combination of multiple zeta values. 

\smallskip

We recall briefly (at the beginning of \S \ref{parintSec} below) 
the standard derivation of the parametric form of
Feynman integrals and its relation to the graph hypersurfaces, 
(see \cite{BjDr}, or the overview given in \cite{Mar}), but first
we recall here the main resulting algebro-geometric setting. 

\medskip
\subsection{Graph hypersurfaces and the polynomial countability question}

Let $\Gamma$ be a finite graph, occurring as a Feynman diagram
for a perturbative scalar massless quantum field theory. We can assume that
$\Gamma$ is connected and one-particle-irreducible (1PI), that is, it cannot be
disconnected by removal of a single edge.
The Kirchhoff polynomial  of $\Gamma$ is defined as
\begin{equation}\label{PsiG}
\Psi_\Gamma = \sum_{T \subseteq \Gamma} \prod_{e \notin T} t_e ,
\end{equation}
where $T \subset \Gamma$ runs through all the spanning trees of $\Gamma$ 
and $t_e$ is an indeterminate corresponding to the edge $e$.  The definition
can be extended to non-necessarily connected graphs, by taking a sum over
spanning forests made of a union of a spanning tree in each connected component.
This a homogenous polynomial in $n=\# E(\Gamma)$ variables, of degree equal to
$b_1(\Gamma)$, hence it defines a projective hypersurface $X_\Gamma$ in 
$\mathbb{P}^{{n-1}}$ given by the set of zeros of $\Psi_\Gamma(t)$. These hypersurfaces
are defined as varieties over $\mathbb{Z}$. Thus, it makes sense to consider the
reductions $X_{\Gamma,p}$ modulo primes, and the counting
of points over finite fields $\F_q$, $q=p^m$, for these varieties. 

\smallskip

The ``polynomial countability" condition for a variety $X$ defined over $\Z$ is the property that
the counting function $N_q(X)= \# X_q(\F_q)$ is a polynomial in $q$. If the motive of the
variety $X$ is mixed Tate, then the variety is polynomially countable (for all but 
finitely many primes) and assuming Tate's conjecture the converse would also hold.
In between the motive and the counting function, one has a universal Euler
characteristic, which is the class $[X]$ in the Grothendieck ring of varieties $K_0(\cV_\Z)$.
Polynomial countability is implied by this class being a polynomial in the Lefschetz
motive $\bL =[\A^1]$, the class of the affine line, which is in turn implied by the motive
of $X$ being mixed Tate.

\smallskip

It was conjectured by Kontsevich  that the graph hypersurfaces $X_\Gamma$ would
always be polynomially countable. This was originally verified for graphs with up to
twelve edges in \cite{Stem}, but  was later proven false in general by Belkale and 
Brosnan \cite{BeBro}, using a powerful universality results for matroids with
which they showed that the classes $[X_\Gamma]$ span a localization of the
Grothendieck ring of varieties, hence they are generally not polynomially countable.

\smallskip

Although the result of Belkale and Brosnan \cite{BeBro} showed that
the graph hypersurfaces are not in general mixed Tate, one also knows that
there are many significant examples (all sufficiently small graphs and several infinite 
families of graphs) for which the varieties in fact happen to be mixed Tate. In the
cases where the mixed Tate property holds, one can then ask a further question of
whether these varieties admit an $\F_1$-structure. 

\smallskip

In addition to the graph hypersurfaces, there are other algebraic varieties
that are naturally associated to the parametric Feynman integrals and that
were recently studied by Esnault and by Bloch, which, unlike the graph
hypersurfaces, are {\em always} mixed Tate motives. It is then even more
natural to ask for these varieties the question of whether (or when) they
admit an $\F_1$-structure. 

\medskip 
\subsection{Parametric Feynman integrals revisited}\label{parintSec}

We present quickly a variant of the usual derivation of the parametric Feynman
integral (as in \cite{BjDr}) and we show that it leads to
the construction of a class of algebraic varieties associated to
graphs, which map naturally to the graph hypersurfaces. These
varieties $\Lambda_\Gamma$ were recently considered by
Bloch \cite{Bloch} and by Esnault: it is shown in \cite{Bloch}
that they are birational covers of the usual graph hypersurfaces whenever
the latter are irreducible (that is, for
all graphs not obtained by gluing two disjoint graphs at a vertex)
and that, unlike the graph hypersurfaces, the $\Lambda_\Gamma$
are always mixed Tate.

\medskip

Recall that the circuit matrix of an oriented graph $\Gamma$ with
a choice of a basis $\{ \ell_j \}$ of $H^1(\Gamma, \Z)$ is the
$\#E(\Gamma)\times b_1(\Gamma)$-matrix defined by
\begin{equation}\label{circuitmat}
\eta_{e,j} =\left\{ \begin{array}{rl} \pm 1 & e\in \ell_j \text{ with same or opposite orientation} \\
0 & e\notin \ell_j . \end{array}\right.
\end{equation}

\medskip

\begin{thm}\label{ParIntLambda}
For a Feynman graph $\Gamma$, let $\cX_\Gamma$ be the locus 
\begin{equation}\label{cXGamma}
\cX_\Gamma=\{ (a,\beta,x)\,|\, a\in \Sigma_n, \, (a,\beta)\in \Lambda_\Gamma, \, x\in Y_{a,m}^c \},
\end{equation}
where
$\Sigma_n$ is the $n$-simplex, 
\begin{equation}\label{Yam}
Y_{a,m}=\{ x \,|\, x^\tau Q_a(x)+m^2=0 \}
\end{equation}
and 
\begin{equation}\label{LambdaGamma}
 \Lambda_\Gamma = \{ (a,\beta) \,|\, Q_a(\beta) =0 \} , 
\end{equation} 
with 
\begin{equation}\label{Qabeta}
(Q_a(\beta))_j = \sum_{e,r} a_e \eta_{e,r} \eta_{e,j} \beta_r ,
\end{equation}
and with $Y_{a,m}^c$ the hypersurface complement.
Then the (unrenormalized) parametric Feynman integral $U(\Gamma)$
with trivial external momenta and non-zero mass is given by the integral
\begin{equation}\label{ParQm}
U(\Gamma)= \int \int \frac{\delta(1-\sum_e a_e) \delta(Q_a(\beta))} {( m^2 + x^\tau Q_a(x))^n} \, da_1 \cdots da_n \, \omega(x,\beta) ,
\end{equation}
where $\omega(x,\beta)$ is the volume form on the fiber in $\cX_\Gamma$
over a point $a\in \Sigma_n$.
\end{thm}

\medskip

\proof
In the usual computation of the parametric Feynman integral
as given in \cite{BjDr}, one starts with a Feynman integral, for a 
Feynman graph $\Gamma$, of the form specified by the Feynman rules, namely
$$ \int \frac{\delta(\sum_{e\in E_{int}} \epsilon_{e,v} k_e + \sum_{e'\in E_{ext}} \epsilon_{e',v} p_{e'})}{q_1(k_1)\cdots q_n(k_n)} \, d^Dk_1 \cdots d^Dk_n, $$
where the constraint in the delta function is momentum conservation
\begin{equation}\label{momcons}
\sum_{e\in E_{int}} \epsilon_{e,v} k_e + \sum_{e'\in E_{ext}} \epsilon_{e',v} p_{e'} =0
\end{equation}
at vertices, with assigned external momenta $p_{e'}$ and the quadrics
$q_e(k_e)=k_e^2+m^2$, and with $D$ the 
spacetime dimension. The matrix $\epsilon_{e,v}$ is the incidence matrix
of the graph $\Gamma$.

\smallskip

Then one performs a change of 
variables,
\begin{equation}\label{varchange}
 k_e = u_e + \sum_{j=1}^\ell \eta_{e,j} x_j , 
\end{equation} 
where $\eta_{e,j}$ is the circuit matrix of the graph $\Gamma$ and $\ell=b_1(\Gamma)$, 
subject to the constraint
\begin{equation}\label{sum0}
 \sum_e a_e u_e \eta_{e,j}=0 , \ \ \ \forall j=1,\ldots,\ell .
\end{equation} 
The $u_e$ are usually taken to be a {\em fixed choice} of a solution to the
resulting equation
\begin{equation}\label{momcons2}
 \sum_{e\in E_{int}} \epsilon_{e,v} u_e + \sum_{e'\in E_{ext}} \epsilon_{e',v} p_{e'}, 
\end{equation}
which follows from \eqref{momcons} because of the orthogonality relation
\begin{equation}\label{matrixorthog}
\sum_e \epsilon_{e,v} \eta_{e,j} =0, \ \ \  \forall v\in V(\Gamma), \, \forall j =1,\ldots, \ell.
\end{equation}

\smallskip

After applying the Feynman trick and performing this change of variables 
in the internal momenta, the Feynman integral becomes of the form
$$ \int \int \frac{\delta(1-\sum_e a_e) \delta(\sum_{e\in E_{int}} \epsilon_{e,v} u_e + \sum_{e'\in E_{ext}} \epsilon_{e',v} p_{e'}) } {(\sum_e a_e (u_e^2 + m^2) + x^\tau Q_a(x))^n} \, da_1 \cdots da_n \, d^D x_1 \cdots d^D x_\ell . $$
There are no mixed terms in the denominator because of the constraint \eqref{sum0}.

\smallskip

If one looks at this same setting, with the change of variables \eqref{varchange} and
the constraint \eqref{sum0}, but setting the external momenta to zero (while keeping
a non-zero mass), one finds that the $u_e$ satisfy the momentum conservation
$$ \sum_{e\in E_{int}} \epsilon_{e,v} u_e =0, $$
that is, they are in the kernel of the incidence matrix. By the orthogonality relation
\eqref{matrixorthog}, one can then express them in the form
$$ u_e = \sum_{r=1}^\ell \eta_{e,r} \beta_r, $$
where the constraint \eqref{sum0} then becomes of the form
\begin{equation}\label{kerQ}
0= \sum_{e,r} a_e \eta_{e,r} \eta_{e,j} \beta_r =Q_a(\beta) .
\end{equation}

While for many choices of $a$ this equation only has the solution $\beta=0$,
which would give the usual choice $u_e=0$ for zero external momenta, it
makes sense here to consider all possible solutions, namely the locus 
given by the complete intersection 
variety $\Lambda_\Gamma$ defined as in \eqref{LambdaGamma}.
This means that one replaces the usual Feynman amplitude
with one where one integrates over all these possible solutions.

This means considering an amplitude of the form \eqref{ParQm}, where
the $n$-form $\omega(x,\beta)$ is the form $d^D k_1\cdots d^D k_n$,
expressed in terms of the $dx_j$ and the $d\beta_j$, after the change
of variables as above, and where,
in the denominator, one has $\sum_e a_e m^2=m^2$, because of the constraint
that the $a_e$ lie on the simplex, and $\sum_e a_e u_e^2 =0$, because
of having set the external momenta $p_{e'}=0$. 
\endproof

\medskip
\subsection{Varieties associated to Feynman integrals in momentum space}

The variant of the parametric Feynman integral described in
Theorem \ref{ParIntLambda} above shows that there are
other interesting algebraic varieties, besides the graph hypersurfaces
$X_\Gamma$, naturally associated to Feynman graphs. We will
focus here especially on the first two cases listed below.

\smallskip
\subsubsection{The complete intersection $\Lambda_\Gamma$}

This is the variety defined by \eqref{LambdaGamma}. It was recently 
introduced and studied by Esnault and by Bloch, in relation to Hodge
structures. It is known, \cite{Bloch}, that $\Lambda_\Gamma$ is always
mixed Tate. Moreover, by writing the Kirchhoff polynomial $\Psi_\Gamma(a)
=\det (Q_a)$, it is shown in \cite{Bloch} that the variety 
$\Lambda_\Gamma \subset X_\Gamma \times \P^{b_1(\Gamma)-1}$, with a
projection $\Lambda_\Gamma \to X_\Gamma$ that is a birational map, when
$X_\Gamma$ is irreducible (when $\Gamma$ is not a union of two disjoint 
graphs glued together at a vertex). 

\smallskip
\subsubsection{The hyperplane arrangement}
The construction of the variety $\Lambda_\Gamma$ is closely related to 
a hyperplane arrangement $\cA_\Gamma$ associated to the graph $\Gamma$, given
by the collection of hyperplanes
\begin{equation}\label{arrHe}
H_e = \{ \beta \, |\, \eta_{e,j} \beta_j =0 \},
\end{equation}
in the rational vector space $\Q^{b_1(\Gamma)}=H^1(\Gamma,\Q)$.

\smallskip
\subsubsection{The hypersurfaces $Y_{a,m}$}
These are the varieties defined by \eqref{Yam}. The presence of the
non-trivial mass parameter $m\neq 0$ means that one can think of
the $Y_{a,m}$ as deformations of the varieties defined by 
the equation $x^\tau Q_a(x)=0$.

\medskip
\subsection{Varieties associated to Feynman integrals in configuration space}

Another interesting algebro-geometric setting for perturbative quantum field
theory arises when one considers Feynman integral calculations in configuration
space rather than momentum space. We refer the reader to the two recent
papers \cite{CeyMa1}, \cite{CeyMa2} for a detailed discussion of this setting
and for the main results. We simply recall here what are the algebraic varieties
that arise in this context.

One starts with a smooth projective variety $X$, which is usually assumed
to be $X=\P^D$, with $D$ the spacetime dimension, but that can be taken
more general. Given a Feynman graph $\Gamma$, inside the product
$X^{\# V(\Gamma)}$, with $V(\Gamma)$ the set of vertices of $\Gamma$,
one considers the diagonals $\Delta_e =\{ x_{s(e)}-x_{t(e)}=0 \}$ and
one defines the graph configuration space to be
\begin{equation}\label{ConfGamma}
{\rm Conf}_\Gamma(X) = X^{\# V(\Gamma)} \smallsetminus \cup_e \Delta_e .
\end{equation}
This variety has a wonderful compactification $\overline{{\rm Conf}}_\Gamma(X)$, 
in the sense of  \cite{DP} and \cite{li}, \cite{li2}, whose geometric and motivic
properties were studied in detail in \cite{CeyMa1}. It is obtained from $X^{\# V(\Gamma)}$
by an iterated sequence of blowups along diagonals. This gives
a completely explicit description of the motive, \cite{li}, \cite{CeyMa1}. 
The varieties $\overline{{\rm Conf}}_\Gamma(X)$ include the Fulton--MacPherson
compactifications, as the case of the complete graph.

\medskip
\section{Constraints from Euler characteristic and Grothendieck class}\label{EulGrSec}

In this section, we start with a very naive point of view on
$\F_1$-geometry, based on the original intuition based on
the existence of a good limiting behavior \eqref{Nq1} for
the counting of the number of points over finite field, when
one lets $q\to 1$.

By taking this simple point of view, one can derive
necessary (but generally far from sufficient) conditions
for the existence of an $\F_1$-structure. 
The simplest such necessary condition can be expressed in terms
of constraints on the Euler characteristic.

We begin by recalling some general facts about the Grothendieck ring
of varieties and the Euler characteristic. Then we discuss some specific
cases of graph hypersurfaces $X_\Gamma$
(or their hypersurface complements $Y_\Gamma$), and of 
varieties $\Lambda_\Gamma$ and $\cA_\Gamma$, associated to
Feynman graphs.

\subsection{The Grothendieck ring of varieties}

We will consider here the case where $R$ is either $\Z$ or $\Q$ of
$\C$ or $\F_q$. We denote by $\cV_R$ the category of quasi-projective
varieties defined over $R$ and by $K_0(\cV_R)$ the Grothendieck
ring of varieties.  The latter is constructed by first considering 
the abelian group of $\Z$-linear combinations  
of isomorphism classes $[X]$ of varieties $X\in \cV_R$ and then
moding out by the subgroup generated by all elements of the form 
$[X] - [Y] - [X\smallsetminus Y]$, where $Y \subseteq X$ is a closed 
subvariety of $X$.  This is known as the scissor -congruence 
or inclusion-exclusion relations. Moreover, $K_0(\cV_R)$ 
is made into a ring by defining the product of two classes 
$[X]$ and $[Y]$ to be the class of the product of $X$ and $Y$ over ${\rm Spec}(R)$, 
that is $[X][Y] = [X \times Y]$, and then extending it by linearity.
The classes $[X]$ in the Grothendieck ring are also referred to as
``virtual motives". The Lefschetz motive is the class $\bL=[\A^1]$ and
one refers to the subring $\Z[\bL]\subset K_0(\cV_R)$ as the
subring of virtual Tate motives. Thus, we say that a variety $X$ is a virtual
Tate motive if its class $[X]=P(\bL)$ for a polynomial $P$ with $\Z$
coefficients, that is, if $[X]\in \Z[\bL]$. 

\subsection{Additive invariants and the Euler characteristic}

The importance of the Grothendieck ring becomes apparent when studying the properties of 
``additive invariants" of algebraic varieties, the prototype example of which is the Euler
characteristic. An additive invariant is a map $\chi: \cV_R \to S$, with $S$ a commutative
ring, with the following properties: 
\begin{enumerate}
\item Isomorphism invariance: $\chi(X) = \chi(Y)$ for all $X \cong Y$.
\item Inclusion-exclusion:  $\chi(X) = \chi(Y) + \chi(X\setminus Y)$, whenever 
$Y \subseteq X$ is a closed subvariety.
\item Multiplicative: $\chi(X \times Y) = \chi(X) \chi(Y)$.
\end{enumerate}

Examples of additive invariants are the topological Euler characteristic when $R=\Q$ or
$\C$, the Deligne--Hodge polynomial when $R=\C$, or the counting of
points over finite fields when $R=\Z$ or $\F_q$. 

Note that the standard terminology ``additive" refers to the behavior with respect to
inclusion--exclusion, but in fact these invariants are also required to be multiplicative
with respect to products, so that they are compatible with the ring structure and not
just with the group structure in the Grothendieck ring.

The Grothendieck ring is in fact universal with respect to additive invariants. 
That is, for any additive invariant $\chi: \cV_R \to S$, there exists a unique ring 
homomorphism $\chi_* : K_0(\cV_R) \to S$ such that the following 
diagram commutes:
$$
\xymatrix{ \cV_R \ar[r]^{\pi} \ar[dr]^{\chi} & K_0(\cV_R) \ar @{.>}[d]^{\chi_*} \\ & S }
$$
where $\pi$ is the map $X \mapsto [X]$. Thus, one can think of the Grothendieck
class $[X]$ as a ``universal Euler characteristic" for the variety $X$.

\medskip
\subsection{Euler characteristic and $\F_1$-points}

The following observation links the Euler characteristic to a necessary
condition for the existence of an $\F_1$-structure over a variety defined
over $\Z$ (see \cite{Deit1}, \cite{Soule}).

\begin{prop}\label{F1chi}
Let $X \in \cV_\Z$ be a variety over $\Z$ such that $[X]\in K_0(\cV_\Z)$ is a 
virtual Tate motive. Then 
\begin{equation}\label{chiNq}
\chi(X_\Q) = \lim_{q\to 1} N_q(X),
\end{equation} 
where $X_\Q$ is a model over $\Q$ of $X$. Thus, a necessary condition
for the existence of an $\F_1$-structure over $X$ is that $\chi(X_\Q)\geq 0$.
\end{prop}

\proof The assumption that $[X]\in K_0(\cV_\Z)$ is a virtual Tate motive 
means that there exists a polynomial $P(t)\in \Z[t]$ such that $[X]=P(\bL)$. 
This implies that, for all but finitely many primes, the counting function satisfies 
$N_q(X)=P(q)$, since the counting is an additive invariant with 
$N_q(\A^1)=N_q(\bL)=q$. 
On the other hand, the Grothendieck ring $K_0(\cV_\Z)$
maps to the Grothendieck ring $K_0(\cV_\Q)$ with $\bL$ mapping to $\bL$,
so that we still have $[X]=P(\bL)$. The Euler characteristic satisfies $\chi(\bL)=1$
and is a ring homomorphism $\chi: K_0(\cV_\Q)\to \Z$, so that we obtain
$\chi(X)=P(1)$, which is the limit of the values $P(q)$ as $q\to 1$. Since the
limit of $N_q(X)$ as $q\to 1$ is interpreted as the counting of points over $\F_1$,
a necessary condition for the existence of an $\F_1$-structure is that this
counting is non-negative.
\endproof

\smallskip

Regardless of the particular flavor of $\F_1$-geometry
one chooses to work with \cite{LL2}, the existence of an $\F_1$-structure on a variety
defined over $\Z$ implies that the motive of the variety $X$ is mixed Tate.
In the following we will restrict our attention to the formulation of  $\F_1$-structures
in terms of torifications given in \cite{LL} and we will show how that implies that
$[X]\in K_0(\cV_\Z)$ is a virtual Tate motive, with a condition on the coefficients
(see Lemma \ref{GrTor} below), which implies the necessary condition 
that $\chi(X_\Q)\geq 0$ as above.

\smallskip
\subsection{Euler characteristic of graph hypersurfaces}

It is often convenient to switch between the projective hypersurfaces $X_\Gamma\subset\P^{n-1}$
and their affine cones $\hat X_\Gamma \subset \A^n$, namely the affine hypersurfaces $\hat X_\Gamma$ defined by the vanishing of the same polynomial $\Psi_\Gamma$ in affine
space $\A^n$. 

\begin{lem}\label{hatXTate}
The class $[X_\Gamma]$ is a virtual Tate motive if and only if the class $[\hat X_\Gamma]$
is a polynomial in $\bL$ divisible by $(\bL-1)$.
\end{lem}

\proof It is easy to check (see \cite{AluMa1}) that the classes in the Grothendieck
ring are related by the simple relation, 
\begin{equation}\label{XhatX}
[X_\Gamma] = \frac{[\hat X_\Gamma] -1}{\bL -1}.
\end{equation}
The statement then follows immediately.
\endproof

\smallskip

We say that an edge $e$ in the graph $\Gamma$ is a bridge if the removal of
$e$ disconnects $\Gamma$ (in particular $\Gamma$ is not a 1PI graph) and
we say that $e$ is a looping egde if the two endpoints of $e$ coincide with the
same vertex of $\Gamma$.

\smallskip

It was shown in \cite{AluMa3} (see also \cite{Stem} for a formulation in terms
of counting functions) that the classes of the affine graph hypersurface
complements satisfy a deletion-contraction relation of the following form:
\begin{equation}\label{delcon}
\begin{array}{ll}
[\A^n \smallsetminus \hat X_\Gamma] = \bL \cdot [\A^{n-1} \smallsetminus \hat X_{\Gamma\smallsetminus e}] =\bL\cdot [\A^{n-1} \smallsetminus \hat X_{\Gamma/e}], & e \text{ bridge} \\[3mm]
[\A^n \smallsetminus \hat X_\Gamma] = (\bL-1) [\A^{n-1} \smallsetminus \hat X_{\Gamma\smallsetminus e}]  = (\bL-1) [\A^{n-1} \smallsetminus \hat X_{\Gamma/e}] & e \text{ looping edge}
\\[3mm]
[\A^n \smallsetminus \hat X_\Gamma] = \bL \cdot [\A^{n-1} \smallsetminus (\hat X_{\Gamma\smallsetminus e} \cap \hat X_{\Gamma/e})] - [\A^{n-1} \smallsetminus \hat X_{\Gamma\smallsetminus e}], & \text{otherwise} 
\end{array}
\end{equation}

\smallskip

\begin{prop}\label{chin}
If $\Gamma$ is not a forest and has at least one bridge 
or one looping edge, then $\chi(X_\Gamma) = n$.
\end{prop}

\proof If $e$ is a bridge in $\Gamma$, and $\Gamma$ is not a forest, 
then  (as observed in \cite{AluMa1} and \cite{AluMa2}) the hypersurface 
complement $\P^{n-2} \smallsetminus X_{\Gamma\smallsetminus e}$ 
is a $\bG_m$-bundle over the product 
$(\P^{n_1-1}\smallsetminus X_{\Gamma_1}) \times 
(\P^{n_2-1}\smallsetminus X_{\Gamma_2})$, 
where $\Gamma\smallsetminus e =
\Gamma_1 \amalg \Gamma_2$ with $n_i =\# E(\Gamma_i)$ for $i=1,2$.
Thus, $[Y_{\Gamma\smallsetminus e}]=(\bL-1) [Y_{\Gamma_1}] \, [Y_{\Gamma_2}]$,
since in the Grothendieck ring the class of the $\bG_m$-bundle is equal to the 
class of the product of base and fiber. Moreover, from the deletion-contraction
relation \eqref{delcon} for the bridge case and the formula  \eqref{XhatX} 
relating the classes of the affine and projective hypersurfaces we obtain
$$ [X_\Gamma ] = \bL [X_{\Gamma\smallsetminus e}] +1, $$
from which we obtain
$$ [Y_\Gamma]= \sum_{k=0}^{n-1} \bL^k -[X_\Gamma] =
\sum_{k=0}^{n-1} \bL^k -1 -\bL [X_{\Gamma\smallsetminus e}] =\bL (
\sum_{k=0}^{n-2} \bL^k -[X_{\Gamma\smallsetminus e}])=\bL [Y_{\Gamma\smallsetminus e}]. $$
Thus, we have $[Y_\Gamma]=\bL (\bL-1) [Y_{\Gamma_1}]\, [Y_{\Gamma_2}]$, which implies,
from the fact that the Euler characteristic is a ring homomorphism and that $\chi(\bL)=1$,
that $\chi(Y_\Gamma)=0$, hence $\chi(X_\Gamma)=\chi(\P^{n-1})=n$.

If $e$ is a looping edge in the graph, then we can rewrite the 
deletion-contraction formula \eqref{delcon} in the equivalent form
\begin{equation}\label{classlooping}
[\hat X_\Gamma]= \bL^{n-1} + (\bL-1) [\hat X_{\Gamma/e}],
\end{equation}
just by writing out the above as $\bL^n -[\hat X_\Gamma]=\bL^n - \bL^{n-1} 
- \bL [\hat X_{\Gamma/e}] + [\hat X_{\Gamma/e}]$.  Using then the
relation \eqref{XhatX} between the classes of the affine and projective
hypersurfaces, we can rewrite this as
\begin{equation}\label{classloopingproj}
[X_\Gamma] = (\bL-1) [X_{\Gamma/e}] +1 + \frac{\bL^{n-1} -1}{\bL-1}.
\end{equation}
We then use the fact that the Euler characteristic is a ring homomorphism
and that $(\bL^{n-1} -1)/(\bL-1)=1+\bL +\cdots + \bL^{n-2}$, so that we
obtain $\chi(X_\Gamma)=(\chi(\bL)-1)\chi(X_{\Gamma/e}) + 1 + \sum_{k=0}^{n-2} \chi(\bL^k)$.
Since $\chi(\bL)=1$, we obtain $\chi(X_\Gamma)=1+(n-1)=n$.
\endproof

\smallskip

This explains the frequent occurrence of the value $\chi(Y_\Gamma)=0$
among sufficienly small graphs, as many have either bridges or looping edges.
It also shows that all these graphs satisfy the necessary condition of Proposition
\ref{F1chi} on the positivity of the Euler characteristic. For graphs that have no
bridges and no looping edges, as observed in \cite{AluMa3}, the deletion-contraction
relation only provides a constraint on the Euler characteristic of the form
\begin{equation}\label{genEulchar}
\chi(X_\Gamma) = n + \chi(X_{\Gamma/e}\cap X_{\Gamma\smallsetminus e}) -
\chi(X_{\Gamma\smallsetminus e}),
\end{equation}
if the graph has $b_2(\Gamma)\geq 2$, and $\chi(X_\Gamma)=n-1$ 
when $b_1(\Gamma)=1$. The latter is the hyperplane case, which is an example
giving an occurrence of the value $+1$ for $\chi(Y_\Gamma)$.

\smallskip

For example, the Euler characteristic was computed explicitly in \cite{AluMa1}
in the case of the banana graphs (graphs with two vertices and $n$ parallel
edges between them), where one has, for $n\geq 3$,
$$ \chi(X_{\Gamma_n}) = n + (-1)^n, $$
and, by \eqref{genEulchar}, one finds $\chi(X_{\Gamma_n / e}\cap 
X_{\Gamma_n\smallsetminus e})=n+(-1)^n -n + n-1+(-1)^{n-1} =n-1$.

\smallskip

\begin{rem}\label{banana}{\rm 
The banana graphs $\Gamma_n$ with $n$ even and $n\geq 3$ give 
examples where $\chi(Y_{\Gamma_n})=-1$.
Thus, in this case, while it is known from \cite{AluMa1} that $Y_{\Gamma_n}=\P^{n-1}\smallsetminus X_{\Gamma_n}$ is a virtual Tate motive, it violates the Euler 
characteristic condition of Proposition \ref{F1chi} for the existence of an $\F_1$-structure.
The hypersurface $X_{\Gamma_n}$ instead is a virtual Tate motive that satisfies the condition 
of Proposition \ref{F1chi}.}
\end{rem}

\smallskip

In all the cases listed above (see Proposition \ref{chin}) where $\chi(X_\Gamma)=n$,
the projective hypersurface complement $Y_\Gamma=\P^{n-1}\smallsetminus X_\Gamma$
satisfies $\chi(Y_\Gamma)=0$. In such cases, assuming that the varieties are virtual Tate
motives, the ``number of points over $\F_1$" is given by the limit
\begin{equation}\label{NF1rG}
N_1(Y_\Gamma) = \lim_{q\to 1} \frac{P_{Y_\Gamma}(q)}{(q-1)^{r_\Gamma}}, 
\end{equation}
where $[Y_\Gamma]=P_{Y_\Gamma}(\bL) \in K_0(\cV_\Z)$, for $P_{Y_\Gamma}(t)\in \Z[t]$ 
and $r_\Gamma$ is the order of vanishing of $P_{Y_\Gamma}(q)$ when $q\to 1$.

\smallskip

\begin{cor}\label{rGloopbridge}
If $\Gamma$ is not a forest and it has either a bridge or a looping edge, then
the order of vanishing $r_\Gamma$ in \eqref{NF1rG} is equal to $r_\Gamma =
1+ s_\Gamma$ where $s_\Gamma$ is as follows.
\begin{enumerate}
\item If $e$ is a bridge, then $(\bL-1)^{s_\Gamma}$ is the maximal power of $\bL-1$
that divides the product $[Y_{\Gamma_1}]\,[Y_{\Gamma_2}]$, with $\Gamma\smallsetminus e=\Gamma_1\amalg \Gamma_2$.
\item If $e$ is a looping edge, then then $(\bL-1)^{s_\Gamma}$ is the maximal power of $\bL-1$
that divides $[Y_{\Gamma/e}]$.
\end{enumerate}
\end{cor}

\proof (1) If $e$ is a bridge, then we use, as in Proposition \ref{chin} the fact that
$[Y_\Gamma]=\bL[Y_{\Gamma\smallsetminus e}]= \bL (\bL-1) [Y_{\Gamma_1}]\, [Y_{\Gamma_2}]$, which gives the result
for $r_\Gamma$. (2) Similarly,  if $e$ is a looping edge, then we have the 
expression \eqref{classloopingproj}
for the class $[X_\Gamma]$ as in Proposition \ref{chin}, which gives
$[X_\Gamma]=(\bL-1)[X_{\Gamma/e}]+1 +\sum_{k=0}^{n-2} \bL^k$. We then
get 
$$ [Y_\Gamma]= \sum_{k=0}^{n-1} \bL^k -[X_\Gamma] =
\bL^{n-1} -1 -(\bL-1) [X_{\Gamma/e}]  $$ $$ = (\bL-1) (\sum_{k=0}^{n-2} \bL^k -  [X_{\Gamma/e}]) =
 (\bL-1)  [Y_{\Gamma/e}]. $$
 So we again obtain the value of $r_\Gamma$ as stated.
 \endproof
 
\medskip

The very simple instances considered here above all lead to values of 
$\chi(Y_\Gamma)$ that are $0$, $+1$ or $-1$. In fact, this is conjectured
to be the always the case.

\begin{conj}\label{chi0pm1conj} {\rm (Aluffi)}
The graph hypersurface complements $Y_\Gamma$ always have
Euler characteristic $\chi(Y_\Gamma)$ which is either $0$ or $1$ or
$-1$.
\end{conj}

This conjecture was so far confirmed by computer calculations on a significant 
number of sufficiently small graphs \cite{Stry}, and for some other classes 
of graphs, using techniques from \cite{Alu5}.
If this conjecture holds, then the Euler characteristics of the graph
hypersurfaces $X_\Gamma$ themselves would always be non-negative,
so this simple necessary condition would always be satisfied. 
  
 \medskip
 \subsection{Constraints from the Grothendieck class}
 
 We now look at the existence of $\F_1$-structures on $\Z$-varieties 
 from a more sophisticated point of view, which is the one developed 
 in the work of Javier L\'opez Pe\~{n}a and Oliver Lorscheid \cite{LL},
 based on the existence of affine torifications. 
 We show that this imposes a very simple necessary condition
 on the form of the class in the Grothendieck ring, which we then discuss
 in some cases of varieties of Feynman graphs.
 
 \smallskip

 \begin{defn}\label{tordef} {\rm (L\'opez Pe\~{n}a and Lorscheid  \cite{LL})}
 A torified scheme is a scheme $X$
 together with a disjoint union of tori $T=\amalg_j \bG_m^{d_j}$ and
 a morphism of schemes $e_X : T \to X$, such that the restriction
 $e_X|_{\bG_m^{d_j}}$ to each torus is an embedding and the morphism
 induces bijections $e_X(k): T(k) \to X(k)$, for every field $k$.
 \end{defn}
 
 \smallskip
 
 \begin{lem}\label{GrTor}
 Let $X$ be a smooth quasi-projective variety defined over $\Z$, which 
 admits a torification $e_X: T \to X$ as above. 
 Then the class $[X]$ in the Grothendieck ring can be written as a
 polynomial $[X] = \sum_{k\geq 0} a_k \bT^k$, in the class $\bT=[\bG_m]=\bL-1$, with
 non-negative integer coefficients $a_k\geq 0$.
\end{lem}

\proof The torification condition implies that the counting function \eqref{NXq} satisfies
$N_X(q) = \sum_j N_{\bG_m^{d_j}}(q)$. The Tate conjecture predicts that the numbers 
$N_X(q)$ determine the motive of $X$, hence its class in the Grothendieck ring, which
should therefore be of the form $[X]=\sum_j [\bG_m^{d_j}]=\sum_j \bT^{d_j}$. 
\endproof

Thus, the form of the class in the Grothendieck ring gives a second
necessary condition for the existence of an $\F_1$-structure, more
refined than the condition on the Euler characteristic. Notice that, if the
class of $X$ in the Grothendieck ring has an expression of the form
$[X]=\sum_k a_k \bT^k$, then the Euler characteristic $\chi(X)=a_0$,
since all the positive dimensional tori have trivial Euler characteristic.

\medskip
\subsection{The case of graph hypersurfaces}

We show in some examples the condition of Lemma \ref{GrTor}
for graph hypersurfaces. As the first example, we consider the
banana graphs, for which the class $[X_\Gamma]$ in the Grothendieck 
ring was computed explicitly in \cite{AluMa1}.

\begin{lem}\label{bananaGrcond}
Let $\Gamma_n$ be the banana graph with two vertices and $n$ parallel
edges between them. The class $[X_{\Gamma_n}]$ satisfies the condition
of Lemma \ref{GrTor}, while the class of the complement $[Y_{\Gamma_n}]=
[\P^{n-1}\smallsetminus X_{\Gamma_n}]$ does not.
\end{lem}

\proof
The class in the Grothendieck ring is a function of
$\bT$ of the form (Theorem 3.10 of \cite{AluMa1}):
\begin{equation}\label{bananaGr}
[X_{\Gamma_n}]= \frac{(1+\bT)^n-1}{\bT} - \frac{\bT^n-(-1)^n}{\bT+1} -n \bT^{n-2}.
\end{equation}
While it appears at first that this expression may contain negative coefficients
of some powers of $\bT$, one can check that in fact it does satisfy the necessary
condition of Lemma \ref{GrTor}. For example, for $n=15$, we find
$$ \begin{array}{rl}
[X_{\Gamma_{15}}]= &
14 + 106\, \bT + 454\, \bT^2 + 1366\, \bT^3 + 3002\, \bT^4 + 5006\, \bT^5 + 6434\, \bT^6 \\
+ &  6436\, \bT^7 + 5004\, bT^8 + 3004\, \bT^9 + 1364\, \bT^{10} + 456\, \bT^{11} + 
 104\, \bT^{12} + \bT^{13}. \end{array} $$
In fact, as already observed in \cite{AluMa1}, we can write the first term in \eqref{bananaGr} as
\begin{equation}\label{PnT}
 [\P^{n-1}] =\frac{\bL^n-1}{\bL-1}= \frac{(1+\bT)^n-1}{\bT} = \sum_{k=1}^n \binom{n}{k}\, 
\bT^{k-1} , 
\end{equation}
while the second term is of the form 
\begin{equation}\label{Tksigns}
 - (\frac{\bT^n-(-1)^n}{\bT+1} + n \bT^{n-2}) =-( \bT^{n-1} +(n-1) \bT^{n-2} +\bT^{n-3} -\bT^{n-4}+\cdots + \pm 1). 
\end{equation} 
Notice then that the coefficients of \eqref{PnT} are always greater than or equal than
the coefficients of the corresponding powers of $\bT$ being subtracted in \eqref{Tksigns},
since the coefficient of $\bT^r$ in \eqref{PnT} greater than or equal to one, in all cases, 
and to $n-1$ in the case of $\bT^{n-2}$.
The class of the hypersurface complement, on the other hand, is given by 
(Corollary 3.13 of \cite{AluMa1})
$$  [Y_{\Gamma_n}]= \bT^{n-1} +(n-1) \bT^{n-2} +\bT^{n-3} -\bT^{n-4}+\cdots + \pm 1, $$
which has some negative coefficients.
\endproof

The deletion-contraction relation for classes in the Grothendieck
ring of graph hypersurface complements, proved in \cite{AluMa3} and
recalled here above in \eqref{delcon}, has as consequence certain
recursion relations for the Grothendieck classes for graphs that
are obtained from one another by iterating certain simple operations,
such as edge splitting or edge doubling. These recursions were
derived in \cite{AluMa3} and can be used to obtain examples
of graphs for which the class $[\hat Y_\Gamma]$ of the affine graph
hypersurface complement satisfies the condition of Lemma \ref{GrTor}.
For example, Proposition 5.11 of \cite{AluMa3} shows that the {\em lemon
wedge graph} $\Lambda_m$ with $m$ sections gives
\begin{equation}\label{melon}
[\hat Y_{\Lambda_m}] =(\bT+1)^{m+1} \sum_{j=0}^m \binom{m-j}{j} \bT^{m-j}.
\end{equation}

\medskip

\begin{ques}\label{F1ques1}
Do graph hypersurfaces $X_\Gamma$, $\hat X_\Gamma$ (or hypersurface complements 
$Y_\Gamma$, $\hat Y_\Gamma$) 
that satisfy the condition of Lemma  \ref{GrTor} admit a torification? Do they admit
an affine torification? Are there natural conditions on the graphs, or
interesting families of graphs for which this is the case? 
\end{ques}

 \medskip
 \subsection{The complete intersection $\Lambda_\Gamma$}
 
 We now discuss briefly the varieties $\Lambda_\Gamma$ associated to
 Feynman graphs. One knows \cite{Bloch} that, unlike the graph hypersurfaces,
 these are mixed Tate for all graphs $\Gamma$. One can then ask whether
 the necessary conditions for the existence of $\F_1$-structures coming from
 the Euler characteristic and the Grothendieck class are satisfied in this case.
 
 It is shown in \cite{Bloch} that the varieties $\Lambda_\Gamma$ have a 
 stratification defined by sets $S_\Gamma^{(m)} \subset \P^{b_1(\Gamma)-1}$
 of the form $S_\Gamma^{(m)}=T_\Gamma^{(m)} \smallsetminus T_\Gamma^{(m+1)}$,
 with each $T_\Gamma^{(m)}$ a union of linear subspaces, which is defined
 as the locus of $\{ \beta \in \P^{b_1(\Gamma)-1}\,|\, \epsilon(\beta)\geq m \}$, where
 $\epsilon(\beta)$ is the codimension of the span of the $\eta_{e,\cdot}$ with
 $\sum_i \eta_{e,i}\beta_i\neq 0$ in $\Q^{b_1(\Gamma)}$. The restriction
 $\Lambda_\Gamma|_{S_\Gamma^{(m)}}$ is a projective bundle over $S_\Gamma^{(m)}$
 with fiber $\P^{n-b_1(\Gamma)+m-1}$. Thus, the class in the Grothendieck ring
 decomposes as
 $$ [\Lambda_\Gamma]=\sum_m [S_\Gamma^{(m)}] \, [\P^{n-b_1(\Gamma)+m-1}]. $$
 Thus, to show that $\Lambda_\Gamma$ satisfies the two necessary conditions, it
 suffices to know that
 $[S_\Gamma^{(m)}]=P_{\Gamma,m}(\bT)$, for a polynomial $P_{\Gamma,m}$ with
 non-negative coefficients, which also implies that $\chi(S_\Gamma^{(m)})\geq 0$ 
 for each stratum.

\begin{ques}\label{F1ques2}
Do the complete intersections $\Lambda_\Gamma$ admit an affine torification? 
\end{ques}

There are explicit cases discussed in \cite{Bloch}, such as the wheel with three
spokes graph, where the variety $\Lambda_\Gamma$ is a projective bundle over
a projective space. This suggests that, at least for certain classes of graphs,
the question above may be checked directly.

 \medskip
 \section{Graph configuration spaces and $\F_1$-geometry}\label{ConfF1Sec}
 
 We now consider the varieties associated to Feynman integrals
 in configuration spaces, as in \cite{CeyMa1} and \cite{CeyMa2}.
 We see that, unlike the case of the varieties of Feynman graphs
 in momentum space, we can say a lot more about $\F_1$-structures
 in the configuration space setting.
 
 \medskip
 \subsection{Euler characteristic and Grothendieck class}
 
 We recall the following notation from \cite{CeyMa1}, see also \cite{li}, \cite{li2}.
 We assume that $\Gamma$ is a finite graph without looping edges (edges whose
 endpoints coincide). 
 
 Recall from  \cite{li2} that a {\it  $\cG_\Gamma$-nest} is 
a collection $\{\gamma_1, \ldots, \gamma_\ell\}$
of biconnected induced subgraphs (induced meaning that the
subgraph has all the same edged between the given subset of
vertices as the original graph), with the property that any two subgraphs
$\gamma$ and $\gamma'$ in the set either have $\gamma \cap \gamma' = \emptyset$,
or $\gamma \cap \gamma' =\{ v \}$, a single vertex, or 
$\gamma \subseteq \gamma'$ or $\gamma' \subseteq \gamma$.
 As in \cite{li2} we use the notation $M_\cN:=\{ (\mu_\gamma)_{\Delta_\gamma \in \cG_\Gamma} \,:\, 1\leq \mu_\gamma \leq r_\gamma -1 , \,\, \mu_\gamma \in \Z \}$ with $r_\gamma= r_{\gamma,\cN} := \dim (\cap_{\gamma'\in \cN: \gamma'\subset \gamma} \Delta_{\gamma'})-\dim \Delta_\gamma$ and
$\| \mu \| := \sum_{\Delta_\gamma \in \cG_\Gamma} \mu_\gamma$.
We also denote by $\Gamma/\delta_{\cN}(\Gamma)$ the quotient
$\Gamma/ / (\gamma_1\cup \cdots \cup \gamma_r)$, where 
$\cN=\{ \gamma_1, \ldots, \gamma_r \}$ is the $\cG$-nest and the quotient double
bar $\Gamma / / \gamma$ means the graph obtained from $\Gamma$ by shrinking
each connected component of the subgraph $\gamma$ to a (different) vertex.

 The class in the Grothendieck ring of the wonderful compactifications
 $\overline{{\rm Conf}}_\Gamma(X)$ was computed in \cite{CeyMa1}:
 it is explicitly of the form
\begin{equation}\label{GrConf}
[ \overline{{\rm Conf}}_\Gamma(X) ] = [X]^{\# V(\Gamma)} +
\sum_{\cN\in \cG_\Gamma\text{-nests}} \,\, 
[X]^{\# V(\Gamma/\delta_{\cN}(\Gamma))} \sum_{\mu \in M_\cN} \bL^{\|\mu\|}.
\end{equation}

\medskip
The necessary conditions for the existence of torifications, in terms of
Euler characteristic and Grothendieck class, can be checked easily as
follows.
 
 \begin{lem}\label{GrEulConf}
 Suppose that the smooth projective variety $X$ is defined over $\Z$ and admits
 an affine torification, which gives it the structure of an $\F_1$-variety in the
 sense of \cite{LL}. Then the wonderful compactifications $\overline{{\rm Conf}}_\Gamma(X)$
 satisfy the necessary conditions for the existence of an $\F_1$-structure, $\chi( \overline{{\rm Conf}}_\Gamma(X)) \geq 0$ and $[\overline{{\rm Conf}}_\Gamma(X)]=\sum_{k\geq 0} a_k \bT^k$ with
$a_k\geq 0$.
 \end{lem}
 
 \proof Both conditions follow immediately from the formula \eqref{GrConf}. In fact, we see
 that 
 \begin{equation}\label{EulConf}
 \chi( \overline{{\rm Conf}}_\Gamma(X)) = \chi(X)^{\# V(\Gamma)} +
 \sum_{\cN\in \cG_\Gamma\text{-nests}} \chi(X)^{\# V(\Gamma/\delta_{\cN}(\Gamma))} 
 \#  M_\cN,
 \end{equation}
 using $\chi(\bL)=1$ and the fact that $\chi$ is a ring homomorphism from
 the Grothendieck ring of varieties to $\Z$. The assumption on $X$ implies
 that $\chi(X)\geq 0$, hence \eqref{EulConf} also gives a non-negative number.
 For the second condition, the assumption on $X$ implies that $[X]=\sum_{j\geq 0} b_j \bT^j$,
 with all the coefficients $b_j \geq 0$. Moreover, all powers $\bL^n$ can
 be written as a polynomial in $\bT$ with non-negative coefficients, using the simple identity
 $$ \sum_{k=0}^n (\bL-1)^k \binom{n}{k} = \bL^n. $$
 It then follows that the second condition is also satisfied.
 \endproof

 \medskip
 \subsection{Affine torifications}
 
 According to the approach to $\F_1$-geometry developed in \cite{LL},
 the existence of a torification is not in itself sufficient for the existence
 of an $\F_1$-structure on a variety defined over $\Z$. However, is the
 torification is in addition {\em affine} then one obtains from it an 
 $\F_1$-{\em gadget} and an $\F_1$-variety, \cite{LL}. The $\F_1$-structure
 can in turn depend on the choice of the torification.
 
 \begin{defn}\label{afftordef} {\rm (L\'opez Pe\~{n}a and Lorscheid  \cite{LL})}
 The {\em affine} condition on a torification $e_X : T \to X$ consists of the
 requirement that there exists an affine covering $\{ \cU_\alpha \}$ of $X$
 that is compatible with the torification, in the sense that, for each affine open
 set $\cU_\alpha$ in the covering, there is a subcollection of tori $T_\alpha$
 in the torification $T$ of $X$, such that the restriction $e_X|_{T_\alpha}: 
 T_\alpha \to \cU_\alpha$ is also a torification.
 \end{defn}
 
 \medskip
 \subsection{Affine torifications and toric bundles}
 
As a preliminary step towards discussing blowups of $\F_1$-varieties,
we show that torifications extend to equivariant projective bundles over
toric varieties.

A vector bundle $\cE$ over a smooth projective toric variety $X$
is said to be equivariant (or toric) if the torus action on $X$ lifts to
an action on $\cE$ that is linear on the fibers. The bundle $\cE$ is
said to have an equivariant structure. A result of \cite{Kly} shows
that $\cE$ is equivariant if an only if the bundles $\gamma^* \cE$ are
all isomorphic to $\cE$, for all elements $\gamma$ in the torus $T$.
For further characterizations and results on toric bundles see 
\cite{HeMuPa} and \cite{BiMuSa}. Vector bundles that are sums of
line bundles admit an equivariant structure, and so do tangent and
cotangent bundle.

\begin{lem}\label{torProj}
Let $X$ be a smooth projective toric variety and let $\cE$ be an
equivariant vector bundle on $X$. Then the projective
bundle $\P(\cE)$ admits an affine torification. 
\end{lem}
 
\proof The equivariant structure on $\cE$ induces a torus action
on the projectivized bundle $\P(\cE)$, so that the projection map
$\pi: \P(\cE)\to X$ is equivariant. Thus, although $\P(\cE)$ is in
general not itself a toric variety, the decomposition of $\P(\cE)$
into orbits of the torus action gives a torification of $\P(\cE)$.
To see that one in fact has an affine torification, notice that the variety $X$
itself has an affine open covering $\{ \cU_\sigma \}$ compatible with the 
torification, by Proposition 1.13 of \cite{LL}, where the open sets
$\cU_\sigma$ are the ones associated to the cones $\sigma$.
The vector bundle $\cE$ is determined, as an equivariant bundle,
by its restrictions $\cE|_{\cU_\sigma}$ and by the equivariant gluings
over the intersections $\cU_{\sigma_i}\cap \cU_{\sigma_j}$. 
Moreover, as observed in \cite{HeMuPa}, the restrictions $\cE|_{\cU_\sigma}$
decompose as sums of line bundles
$$ \cE|_{\cU_\sigma} \simeq \cL_1 \oplus \cdots \oplus \cL_r. $$
In the case of a sum of line bundles, one knows that the
projectivization $\P(\cE|_{\cU_\sigma})$ is in fact a toric variety,
and this implies that the torification is affine. The equivariant
gluing guarantees the compatibility of the torifications over the
intersections  $\cU_{\sigma_i}\cap \cU_{\sigma_j}$ so that one
obtains an affine torification for $\P(\cE)$.
\endproof 
 
 \medskip
 \subsection{$\F_1$-structures and blowups}
 
 We now consider conditions under which one can perform blowups 
 of $\F_1$-varieties.
 
 \begin{lem}\label{F1blowup}
Let $X$ be a torified smooth projective variety with an affine torification 
$e_X: T_X \to X$. Let $Y \subset X$ be a smooth subvariety which is
toric, with an affine torification $e_X|_{T_Y}: T_Y \to Y$ that consists of 
a subcollection of tori $T_Y$ in the torification $T_X$ and which are
orbits of the torus action on $Y$. Assume moreover that the normal
bundle $\cN_X(Y)$ has an equivariant structure. 
Then the blowup ${\rm Bl}_Y(X)$ of $X$ along $Y$ has an affine torification. 
 \end{lem}

\proof The exceptional divisor $E$ of the blowup can be identified with
the projectivization $\P(\cN_X(Y))$. Thus, by Lemma \ref{torProj} this
admits an affine torification. Moreover, the compatibility between the
affine torification of the toric variety $Y$ and the affine torification of
the (possibly non-toric) ambient variety $X$ imply that the affine torification
obtained in this way on $E$ is compatible with the affine torification of $X$,
so that one obtains an affine torification of  ${\rm Bl}_Y(X)$.
 \endproof

 \medskip
 \subsection{Wonderful compactifications and $\F_1$-structures}\label{wonderSec}
 
We now consider the wonderful compactifications $\overline{\rm Conf}_\Gamma(X)$,
in the case where $X=\P^D$, used in \cite{CeyMa1}, \cite{CeyMa2} to compute
Feynman integrals in configuration space.

\smallskip

We recall briefly the blowup construction of  $\overline{\rm Conf}_\Gamma(X)$,
see  \cite{CeyMa1} and \cite{li}, \cite{li2}.
Let $\cG_{k,\Gamma}\subseteq \cG_\Gamma$ be the subcollection of 
{\em biconnected} induced subgraphs with $\# V(\gamma)=k$.
Let  $Y_0=X^{\# V(\Gamma)}$ and let 
$Y_k$ be the blowup of $Y_{k-1}$ along the (iterated) dominant transform 
of $\Delta_{\gamma} \in \cG_{n-k+1,\Gamma}$.  If $\Gamma$ is itself biconnected, 
then $Y_1$ is the blowup of $Y_0$ along the deepest diagonal $\Delta_\Gamma$, 
otherwise $Y_1=Y_0$. (Similarly, we have $Y_k=Y_{k-1}$ whenever there are 
no biconnected induced subgraphs with exactly $n-k+1$ vertices.)
The resulting sequence of blowups
\begin{equation}\label{Yks}
Y_{n-1} \to \cdots \to Y_2 \to Y_1 \to X^{\#V(\Gamma)},
\end{equation}
where $n=\# V(\Gamma)$, does not depend on the order in which 
the blowups are performed, for $\gamma \in \cG_{n-k+1,\Gamma}$,
for a fixed $k$. The variety $Y_{n-1}$ obtained in this way is   
the wonderful compactification $\overline{\rm Conf}_\Gamma(X)$.

\smallskip

\begin{lem}\label{afftorXV}
Let $X$ be a toric variety and let $\Delta_\gamma =\cap_{e\in E(\gamma)} \Delta_e$
be the diagonal of a biconnected induced subgraph $\gamma\subseteq \Gamma$
in $X^{\# V(\Gamma)}$. Then $\Delta_\gamma$ is toric and there is an affine
torification of $X^{\# V(\Gamma)}$ that induces a compatible affine torification on 
$\Delta_\gamma$.
\end{lem}

\proof The product of toric varieties is toric, hence $X^{\# V(\Gamma)}$ is a
toric variety. Since the graph $\gamma$ is connected, we have (Lemma 1 of \cite{CeyMa1})
$\Delta_\gamma \simeq X^{\# V(\Gamma/\gamma)}$, so this is also a toric variety, and
we can decompose
$X^{\# V(\Gamma)}\cong X^{\# V(\Gamma/\gamma)} \times X^{\# V(\gamma)}$,
with compatible affine torifications. 
\endproof

We then have the following result.

\begin{thm}\label{ConfF1}
Let $X=\P^D$. 
Then, for all Feynman graphs $\Gamma$ the configuration
spaces $\overline{\rm Conf}_\Gamma(X)$ admit an affine torification, hence 
they admit a structure of $\F_1$-variety.
\end{thm}

\proof
Recall that (Definition 2.7 of \cite{li2}) for a blowup $\pi: {\rm Bl}_Z(Y) \to Y$ 
of a smooth subvariety $Z\subset Y$ of a smooth projective variety $Y$, 
the dominant transform of an irreducible subvariety $V$ of $Y$ is the 
proper transform of $V$, if $V$ is not contained in $Z$, and the scheme-theoretic 
inverse image $\pi^{-1}(V)$ if $V\subset Z$. 
Let $Y_k$ be the $k$-th step of the iterated blowup construction \eqref{Yks}
of $\overline{\rm Conf}_\Gamma(X)$ and let $\Delta^{(k)}_\gamma$ be the
dominant transform of the diagonal $\Delta_\gamma$ of a subgraph 
$\gamma \in \cG_{n-k+1,\Gamma}$. We show that, at each step, 
there is an affine torification on $Y_k$ that induces a compatible 
affine torification on the $\Delta^{(k)}_\gamma \subset Y_k$.
To see this, we need to check that,
at each stage of the blowup construction, one can repeatedly apply Lemma \ref{torProj}
and Lemma \ref{F1blowup}. Thus, we need to ensure that, at each step,
the projectivized normal bundles of the blowup loci have an equivariant structure.
An explicit description of these projectivized normal bundles 
was given in \S 2.5 of \cite{CeyMa1}, in terms of {\em screen configurations},
as in \S 1 of \cite{FM}. Namely, one has the following setting. Let $\bT_x=T_x(X)$,
the tangent space of $X$ at $x$. The screen configuration space of a 
subgraph $\gamma$ at $x$ is  $\P(\bT_x^{V(\gamma)}/\bT_x)$. Notice that the
projectivized bundle $\P(\bT^{V(\gamma)}/\bT)$ admits an equivariant structure,
as it is built out of copies of the tangent bundle $\bT=T(X)$, which is an equivariant
vector bundle. The blow up loci in the iterative construction of 
$\overline{\rm Conf}_\Gamma(X)$ are the dominant transforms of the
diagonals $\Delta_\gamma$ in $\cG_{n-k+1,\Gamma}$ and at each stage
the exceptional divisor is identified with the projectivized normal bundle
$$ \P(\cN(\Delta_\gamma \subset \cap_{\gamma'\in \cN\,:\, \gamma'\subsetneq \gamma} \Delta_{\gamma'})), $$
where $\cN$ is a $\cG$-nest containing $\gamma$ (see Proposition 3 of \cite{CeyMa1}), or
with $\P(\cN(\Delta_\gamma\subset X^{V(\Gamma)})$ if 
$\{ \gamma'\in \cN\,:\, \gamma'\subsetneq \gamma\}=\emptyset$.
One also has an identification (Theorem 1 of \cite{CeyMa1})
$$ \P(\cN(\Delta_\gamma \subset \cap_{\gamma'\in \cN\,:\, \gamma'\subsetneq \gamma} \Delta_{\gamma'})) \simeq \P(T(X^{V(\Gamma//(\gamma_1\cup\cdots \cup \gamma_r))})/
T(X^{V(\Gamma/\gamma)})), $$
where $\cN=\{ \gamma_1, \ldots, \gamma_r \}$ and the latter can in turn be identified with
$$ T(X^{V(\Gamma//(\gamma_1\cup\cdots \cup \gamma_r))})/
T(X^{V(\Gamma/\gamma)} \simeq \bT^{V(\gamma//(\gamma_1\cup\cdots \cup \gamma_r))}/\bT ,$$
hence by the previous observation it carries an equivariant structure.
One concludes using Lemma \ref{torProj} and 
Lemma \ref{F1blowup} that at each stage of the construction, the
blowups then have an affine torification.
One can then prove the statement inductively. At the first step $Y_0=X^{\#V(\Gamma)}$
and one has a compatible affine torification on the deepest diagonal $\Delta_\Gamma$
by Lemma \ref{afftorXV}. (We can assume $\Gamma$ itself is biconnected for simplicity.)
Then the general argument here above shows that the dominant transforms 
in $Y_1={\rm Bl}_{\Delta_\Gamma}(X^{\#V(\Gamma)})$ of the
diagonals $\Delta_\gamma$, with $\gamma \in \cG_{n-1,\Gamma}$ has an affine
torification compatible with the affine torification obtained on the blowup $Y_1$.
Similarly, if we assume that the result holds for all the $Y_r$ with $r\leq k-1$, we
can apply the general argument described above to obtain compatible
regular affine torifications on the blowup $Y_k$ and on the dominant transforms
of the  $\Delta_\gamma$, with $\gamma \in \cG_{n-k+1,\Gamma}$.
We can then conclude that the wonderful compactifications $\overline{\rm Conf}_\Gamma(X)$
are $\F_1$-varieties in the sense of \cite{LL}.
\endproof

The result can be extended similarly to cases where $X$ is
a smooth projective toric variety, though the case of $X=\P^D$ is
the one that is physically significant for quantum field theory in
configuration spaces.

\medskip
\subsection{$\F_1$-geometry and the moduli space of curves}

We have discussed in the previous section the existence of
affine torifications on the wonderful compactifications $\overline{\rm Conf}_\Gamma(X)$
of the graph configuration spaces. It is well known that, when $\Gamma$ is the complete
graph $\Gamma_n$, this is just the Fulton--MacPherson space $X[n]$. 
Moreover, it is known that the moduli space $\bar\cM_{0,n+1}$ can be realized as a 
subscheme of $X[n]$ for $X$ any smooth curve, \cite{CheGibKra}. 
Thus, it is natural to ask the question of whether
the method discussed in the previous section can be adapted to show that the
moduli spaces $\bar\cM_{0,n}$ are defined over $\F_1$.

\smallskip

There are various constructions of $\bar\cM_{0,n}$. The one that appears
most directly useful from our point of view is the one described in \cite{CheGibKra},
as part of a larger family of varieties $T_{d,n}$, which are compactifications of
the configuration spaces of $n$ distinct points in $\A^d$ up to translation
and homothety, with the boundary strata parameterizing $n$-pointed stable rooted
trees of $d$-dimensional projective spaces, so that $T_{1,n}=\bar\cM_{0,n}$.

\smallskip

\begin{thm}\label{M0nF1}
The moduli spaces $\bar\cM_{0,n}$ and the varieties
$T_{d,n}$ of \cite{CheGibKra}, for $d\geq 2$ and $n\geq 2$,  are defined over $\F_1$ 
in the sense of \cite{LL}, that is, they have an affine torification. 
\end{thm}

\proof We follow the steps in the iterated blowup construction of \cite{CheGibKra}
and we show that we can apply to them the same results of \S \ref{ConfF1Sec} above.
The varieties $T_{d,n}$ defined in \cite{CheGibKra}, with $T_{1,n}=\bar\cM_{0,n}$
are obtained as an iterated sequence of blowups (Theorem 3.3.1 of \cite{CheGibKra})
\begin{equation}\label{Fdni}
T_{d,n+1}=F_{d,n}^n \to \cdots \to F^1_{d,n}\to F^0_{d,n}=T_{d,n},
\end{equation}
where $F^1_{d,n}\to F^0_{d,n}$ is a projective bundle morphism, and the
other morphisms $F^{i+1}_{d,n} \to F^i_{d,n}$ are blowups along a union
of subvarieties isomorphic to products $T_{d,n-i}\times T_{d,i+1}$. Thus,
$T_{d,n+1}$ is obtained as a sequence of blowups from a projective
bundle over $T_{d,n}$. By Proposition 3.4.1 of \cite{CheGibKra} we also
know that $T_{1,3}\simeq \P^1$, while for $d>1$ we have $T_{d,2}\simeq \P^{d-1}$.
Since $\bar\cM_{0,n}$ is of dimension $n-3$, we are interested in the cases with
$n>3$, and we can take $\bar\cM_{0,4}=T_{1,3}\simeq \P^1$ as the starting point for
the construction. Similarly, we use $T_{d,2}\simeq \P^{d-1}$ as the starting
point for the case of $d>1$. To obtain the result, we again need to check that
at each stage of the iterated blowup construction the normal bundles have an
equivariant structure. To this purpose we can use the fact that the construction
of the spaces $F^i_{d,n}$ can be done inside strata of the compactification of
the Fulton--MacPherson configuration spaces \cite{FM}. In fact, it is shown
in \S 3 of \cite{CheGibKra} that, for a smooth variety $X$ of dimension $d$
(which we can take to be $\P^d$ for convenience), the spaces $T_{d,n}$
are realized as the fibers of the projection map 
$$ X[n]=\overline{\rm Conf}_{\Gamma_n}(X) \to \Delta_{\Gamma_n}\simeq X, $$
over the boundary stratum $D_{\Gamma_n}$ of $X[n]$ that corresponds to the
diagonal $\Delta_{\Gamma_n}$. Therefores, as explained in \S 3 of \cite{CheGibKra}, 
the normal bundles can again be described in terms of screen configurations, 
so we can apply the same argument of Theorem \ref{ConfF1} to show that the 
normal bundles have an equivariant structure so that we can apply
Lemma \ref{torProj} and Lemma \ref{F1blowup}, to obtain affine torifications 
on the $T_{d,n}$ and in particular on $\bar\cM_{0,n}$.
\endproof

\medskip
\section{Hyperplane arrangements and $\F_1$-geometry}\label{HypSec}
 
We discuss in this section a question about $\F_1$-geometry, which
is aimed at identifying other natural and sufficiently interesting
classes of varieties that can, under suitable conditions, carry 
$\F_1$-structures. We identify the hyperplane arrangements
as a sufficiently broad class of varieties (in the specific
sense described in \S \ref{HypMTsec} below), for which
the question of the existence of an $\F_1$-structure (in the sense
of an affine torification as in \cite{LL}) is related to well known
invariants.

\medskip
 \subsection{Hyperplane arrangements and mixed Tate motives}\label{HypMTsec}
 
 One reason why it is especially interesting, from the point of view
 of $\F_1$-geometry, to identify necessary and sufficient conditions
 for hyperplane arrangements to admit an $\F_1$-structure is that
 these varieties are very general from the point of view of mixed
 Tate motives.
 
 \smallskip
 
 In fact, one knows that in all flavor of $\F_1$-geometry, the existence
 of an $\F_1$-structure implies that the motive of the variety is mixed Tate.
 One can therefore formulate the question of characterizing which
 mixed Tate motives over $\Z$ admit an $\F_1$-structure. 
 
 \smallskip
 
 This first requires a good characterization of mixed Tate motives
 in terms of appropriate ``generators". There are several conjectures
 to that extent (see \cite{Brown}), but we focus here especially on the 
 one proposed by Beilinson--Goncharov--Schechtman--Varchenko in
\cite{BGSV}, which predicts that, in the number field case, all the extensions 
${\rm Ext}^1(\Q(0),\Q(n))$  in the category of mixed Tate motives are realized 
in terms of  hyperplane arrangements in general position. This was further
elaborated, more recently, by Madhav Nori \cite{Nori}, who showed that the
construction described in \cite{BGSV} of a graded Hopf algebra of hyperplane 
arrangements  indeed determines a category of mixed Tate motives. 

\smallskip

More precisely, one can formulate the conjecture in the following way 
(\cite{BGSV}, \cite{Brown}).

\begin{conj}\label{HypTate} {\rm (Beilinson--Goncharov--Schechtman--Varchenko)}
For a number field $k$, given a choice of hyperplanes $L_i$, $M_i$, $i=0,\ldots,n$ 
in $\P^n$ in general position, defined over $k$, one obtains a mixed Tate motive
\begin{equation}\label{motHyp}
 \m(\P^n\smallsetminus M, L\smallsetminus (L\cap M)) 
\end{equation} 
whose realizations give the middle dimensional relative cohomology 
$$H^n(\P^n\smallsetminus M, L\smallsetminus (L\cap M)).$$ Then the subcategory
$\cC$ of the category $\cM\cT(k)$ of mixed Tate motives over $k$ generated by
the motives of the form \eqref{motHyp} is in fact all of $\cM\cT(k)$.
\end{conj}

Assuming this conjecture, one can then try to characterize which mixed Tate
motives over $\Z$ descend to $\F_1$, by identifying conditions under which
hyperplane arragements defined over $\Z$ carry an $\F_1$-structure. Notice that,
while the conjecture above is formulated for motives over a number field, 
as explained in \S 6 of \cite{Gon}, the abelian category of 
mixed Tate motives over the ring of integers in a number field can be
obtained by considering a sub-Hopf-algebra of the Hopf algebra that
defines the abelian category of mixed Tate motives over the number field.
  
\medskip
\subsection{Necessary conditions for $\F_1$-structures} 
 
 The necessary condition based on the Grothendieck class (which obviously
 implies the one based on the Euler characteristic) admits a nice interpretation in
 the case of hyperplane arrangements. 
 
 We first recall some standard terminology
 in hyperplane arrangements, \cite{Alu3}, \cite{OrTe}. 
 Given a hyperplane arrangement $\cA=\{ H_j \}_{j=1}^N$ of $N$ hyperplanes in
$\P^n$, the central arrangement $\hat\cA= \{ \hat H_j \}_{j=1}^N$ is the associated
arrangement of hyperplanes in affine space $\A^{n+1}$ through the origin. The characteristic
polynomial $\chi_{\hat\cA}(t)$ of the central arrangement is defined as
\begin{equation}\label{chiAt}
\chi_{\hat\cA}(t) = \sum_{x\in \cL(\hat\cA)} \mu(x)\, t^{\dim(x)},
\end{equation}
where $\cL(\hat\cA)$ is the poset of the $\cap_{j\in J} \hat H_j$, with $J\subseteq \{ 1, \ldots, N\}$,
ordered by reverse inclusion, and $\mu(x)$ is the M\"obius function of $\cL(\hat\cA)$.
 
  \begin{lem}\label{HypGr}
 Let $\cA=\{ H_j \}_{j=1}^N$ be a hyperplane arrangement with $A=\cup_{j=1}^N H_j$.
 Then the Grothendieck class satisfies $[ A ]=\sum_{k\geq 0} a_k \bT^k$, with $a_k \geq 0$, if
 and only if the M\"obius function $\mu(x)$ of $\cL(\hat\cA)$ satisfies, for all $k\geq 1$
 \begin{equation}\label{MobApos}
 \sum_{x:\, \dim(x)\geq k} \mu(x) \, \binom{\dim(x)}{k} \leq \binom{n+1}{k}.
 \end{equation}
 \end{lem} 
 
 \proof  This is a direct consequence of Theorem 1.1 of \cite{Alu3}, according
 to which the class in the Grothendieck ring of the arrangement complement satisfies
 $$ [ \P^n \smallsetminus A ] = \frac{\chi_{\hat\cA}(\bL)}{\bL-1}. $$
 This means that
 $$ [A] = [\P^n]- \frac{\chi_{\hat\cA}(\bL)}{\bL-1}
 =\frac{(1+\bT)^{n+1}-1 -\chi_{\hat\cA}(\bT+1)}{\bT} .
 $$
 Thus, the condition that  $[ A ]=\sum_{k\geq 0} a_k \bT^k$, with $a_k \geq 0$
 can be formulated as the stated condition, using the fact that
 $$ \chi_{\hat\cA}(t+1) =\sum_{x\in \cL(\hat\cA)} \mu(x)\,( t+1)^{\dim(x)} =
 \sum_{x\in \cL(\hat\cA)} \mu(x)\, \sum_{r=0}^{\dim(x)} \binom{\dim(x)}{r} \,  t^r, $$
 so that the above can be written equivalently as
 $$ [A] = \frac{1}{\bT} \sum_{k=1}^{n+1}\left( \binom{n+1}{k} - \sum_{x:\, \dim(x)\geq k} \mu(x)
 \, \binom{\dim(x)}{k} \right) \, \bT^k, $$
 from which \eqref{MobApos} follows.
 \endproof
 
 Notice that the positivity condition described here for the coefficients $a_k$ of the
 class $[A] =\sum_{k\geq 0} a_k \bT^k$ is the same as the condition given in Proposition 6.1
 of \cite{Alu3} for the positivity of the coefficients of the Chern--Schwartz--MacPherson 
 class $c_{SM}(A)$. In fact, the term $x=0$ in the sum on the left-hand-side of \eqref{MobApos}
 has $\mu(0)=1$, hence it cancels the term on the right-hand-side of \eqref{MobApos}, so that
 the same condition can be expressed as the condition that all the coefficients of the
 polynomial $-\sum_{x\neq 0} \mu(x) (t+1)^{\dim(x)}$ are non-negative, which is the
 one given in Proposition 6.1 of \cite{Alu3} for the CSM class. There is an {\em a priori}
 reason for this coincidence, namely the result of Proposition 2.2 of \cite{AluMa1},
 which shows that, for all subvarieties of projective spaces that are obtained as
 unions, intersections, complements, differences of linearly embedded subspaces,
 the CSM class can be obtained from the class in the Grothendieck ring by
 formally replacing each power $\bT^r$ with the class in the Chow group 
 of a linear subspace $\P^r$. We return to discuss the relevance of CSM classes to
 questions about $\F_1$-geometry in \S \ref{ChernSec} below.

\begin{ques}\label{F1ques4}
Is there a combinatorial condition on hyperplane arrangements that is 
{\em sufficient} for the existence of an affine torification on the
arrangement or on its complement?
\end{ques} 
 
\medskip

Coming back to the examples motivated by varieties arising in perturbative
quantum field theory, one can formulate the following more specific question. 
 
\begin{ques}\label{F1ques3}
Do the hyperplane arrangements $\cA_\Gamma$ associated to Feynman graphs 
admit an affine torification? 
Are there natural conditions on the graphs, or interesting families of graphs
for which this is the case?
\end{ques} 

Notice that, in general, one has examples of hyperplane arrangements
that do not satisfy the necessary condition on the Euler characteristic
(Example 6.4 of \cite{Alu3}), so in the case of the arrangements $\cA_\Gamma$
associated to Feynman graphs one seeks conditions on the graph that
might relate to the existence of a torification of the arrangement.

 \medskip
 \subsection{Arrangements, wonderful compactifications, and torifications}
 
 Given a hyperplane arrangement $\cA$, besides the varieties $A$ 
 and $\P^n\smallsetminus A$, there is another class of varieties that
 are interesting to consider, which are the wonderful compactifications,
 in the sense of \cite{DP}, of the complement $\P^n\smallsetminus A$.
 We will denote them here by $\overline{Conf}_{\cA}(\P^n)$, by analogy to
 the notation used above for the graph configuration spaces. Notice that
 the spaces $\overline{Conf}_\Gamma(X)$ describe above do not correspond
 to hyperplane arrangements, because the linear spaces $\Delta_e$
 are complete intersections $\Delta_e=\{ x_{s(e)}=x_{t(e)} \}=\cap_{k=1}^D \{ 
 x_{s(e),k}=x_{t(e),k} \}$, hence they are of codimension $D$ in $X^{\#V(\Gamma)}$.
However, an argument similar to the one used in \S \ref{wonderSec} to 
obtain torifications of  $\overline{Conf}_\Gamma(X)$ can be developed for
the wonderful compactifications $\overline{Conf}_{\cA}(\P^n)$ of \cite{DP},
for hyperplane arrangements, leading to conditions on the existence of
$\F_1$-structures on these varieties.
     
 \medskip
 \section{Algebro-geometric Feynman rules and embedded $\F_1$-structures}\label{ChernSec}
 
In \cite{AluMa2}, a polynomial invariant of the graph hypersurfaces 
$\A^n\smallsetminus \hat X_\Gamma$ is defined using the Chern--Schwartz--MacPherson 
characteristic classes of singular varieties. The polynomial $C_\Gamma(T)$ satisfies the
``algebro-geometric Feynman rule" property of being multiplicative over disjoint union
of graphs, $C_{\Gamma_1\cup \Gamma_2}(T)=C_{\Gamma_1}(T) C_{\Gamma_2}(T)$.
Moreover, it contains the Euler characteristic of the {\em projective} graph hypersurface
complement $\chi(\P^{n-1}\smallsetminus X_\Gamma)$ as the coefficient of the degree one
term. Thus, it can be thought of as a generalization of this Euler characteristic that
satisfies the multiplicative property, which fails for  $\chi(\P^{n-1}\smallsetminus X_\Gamma)$ 
itself. 

\smallskip

We recall here some basic properties of characteristic classes of singular
varieties and the definition of the invariant $C_\Gamma(T)$. We then give
an interpretation of $C_\Gamma(T)$ in terms of Euler characteristic using
a result of Aluffi \cite{Alu6}, and we explain its relevance to 
the $\F_1$-geometry point of view.

\medskip
\subsection{Chern classes of singular varieties}

We recall here briefly a few facts about characteristic classes of singular
varieties that we need to use in the following. For a detailed introduction
to the subject, we refer the reader to \cite{Alu1}.

\smallskip

Let $V$ be a variety over $\bK$. For any closed subvariety $S$ we define $\mbf{1}_S$ to be the function from $V$ to $\Z$ that is $1$ on $S$ and $0$ on $V \setminus S$. Let $\cC:\cV_\bK^p \to Ab$ be the functor from $\bK$-varieties and proper maps to abelian groups where $\cC(V)$ is the abelian group of $\Z$-linear combinations $\sum m_S \mbf{1}_S$ over $S$ closed subvarieties. These are called constructible functions. 
For $f:V \to W$, $\cC(f)(\mbf{1}_S)$ is defined by 
$$
\cC(f)(\mbf{1}_S)(p) = \chi(f^{-1}(p) \cap S),
$$
and extending by linearity. 

\smallskip

Grothendieck and Deligne conjectured the existence of a map $c_*$ from $\cC(V)$ to the homology of $V$ satisfying the following properties for any constructible functions $\alpha$ and $\beta$ on $V$ and any $f:V \to W$:
\begin{enumerate}
\item $c_*(\cC(f)(\alpha)) = f_*c_*(\alpha)$
\item $c_*(\alpha + \beta) = c_*(\alpha) + c_*(\beta)$
\item $c_*(\mbf{1}_V) = c(V) \cap [V]$, when $V$ is smooth,
\end{enumerate}
where $c(V)$ denotes the usual Chern class of $V$ and $c(V) \cap [V]$ is the image of $c(V)$ in homology under Poincare duality. Any such $c_*$, if it exists, is unique. Schwartz and MacPherson independently constructed such a $c_*$, showing existence. (In fact, the classes
defined by M.H.~Schwartz were introduced prior to the functorial formulation, 
and later proved to satisfy the properties above by Brasselet and Schwartz.)
One defines $\csm(V)$ for any (not necessarily smooth) $V$ to be $c_*(\mbf{1}_V)$, see \cite{Alu1} \cite{macpherson}.

\medskip
\subsection{The CSM Feynman rule}

We recall here the definition and properties of the invariant $C_\Gamma(T)$ constructed
in \cite{AluMa2}.

\smallskip

Given a locally closed subset $\hat X$ of $\P^N$, one can write the CSM class in the Chow group $A_*(\P^N)$ as
\begin{equation}\label{CSMakPk}
 c_*(\mbf{1}_{\hat X})=a_0 [\P^0] + \cdots a_N [\P^N] . 
\end{equation}
One then introduces an associated polynomial of the form
\begin{equation}\label{GXT}
G_{\hat X}(T):= a_0 + a_1 T + \cdots + a_N T^N , 
\end{equation}
obtained by formally replacing the class $[\P^k]$ in $A_*(\P^N)$ with the variable $T^k$.
This polynomial can be equivalently written as $G_{\hat X}(T) = \sum_{k\geq 0} c_k T^k$,
where the coefficient $c_k$ is the degree of the $k$-dimensional piece of the CSM class 
of $\hat X$. In particular $G_{\hat X}(0) = \chi(\hat X)$. In the case where $\hat X=\A^N$, 
we have $G_{\A^N}(T)=(T+1)^N$. 

\smallskip

One then defines the polynomial $C_\Gamma(T)$ as the difference
\begin{equation}\label{CGammaT}
C_\Gamma(T) = G_{\A^n}(T) - G_{\hat X_\Gamma}(T),
\end{equation}
for $n=\# E(\Gamma)$. This is a monic polynomial of degree $n$. 
Here the affine graph hypersurface $\hat X_\Gamma \subset \A^n$,
which is the affine cone over the projective $X_\Gamma\subset \P^{n-1}$ is seen as 
embedded in $\P^n$.  It is proved in
\cite{AluMa2} that if $\Gamma = \Gamma_1 \cup \Gamma_2$ is a disjoint union of
graphs, then the polynomial behaves multiplicatively,
\begin{equation}\label{CTprod}
C_\Gamma(T) = C_{\Gamma_1}(T) C_{\Gamma_2}(T) ,
\end{equation}
and it contains the Euler characteristic of the projective hypersurface complement as a coefficient,
\begin{equation}\label{CGTprimechi}
C_\Gamma^\prime(0) = \chi(\P^{n-1}\smallsetminus X_\Gamma).
\end{equation}
The multiplicative property \eqref{CTprod} means that $C: \cH \to \Z[T]$ is 
a {\em ring homomorphism} from the Hopf algebra of graphs to the ring of polynomial.
This is the minimal algebraic requirement needed for a ``Feynman rule", in the setting 
of Hopf algebra-based renormalization, see \cite{AluMa2}.  For a more detailed
discussion of this multiplicative property, see also \cite{Alu6}.

\medskip
\subsection{Chern classes and Euler characteristics}

Let $\hat X\subset \P^N$ be an embedded subvariety and set
\begin{equation}\label{EulcharHi}
\chi_k := \chi(X \cap H_1\cap\cdots\cap H_k),
\end{equation}
where $\chi$ is the topological Euler characteristic and
$H_1, \ldots, H_k$ are general hyperplanes in $\P^N$. One also defines
\begin{equation}\label{chiXt}
\chi_{\hat X} (T) = \sum_{k\geq 0} \chi_k \, T^k.
\end{equation}
The value at zero is the usual topological Euler characteristic, $\chi_{\hat X}(0)=\chi(\hat X)$.

\smallskip

It was proved by Aluffi \cite{Alu6} that the polynomial $G_{\hat X}(T)$
defined in \eqref{GXT} above can be expressed in terms of the polynomial
$\chi_{\hat X}(T)$ of \eqref{chiXt}, and conversely, according to the formulae
\begin{equation}\label{chiandG}
\chi_{\hat X}(T) = \frac{T \,\, G_{\hat X}(T-1) - G_{\hat X}(0)}{T-1} \ \ \  \text{ and } \ \ \
G_{\hat X}(T) = \frac{T \,\, \chi_{\hat X}(T+1) + \chi_{\hat X}(0)}{T+1}.
\end{equation}

\smallskip

This result means that the Chern--Schwartz--MacPherson class can be
computed completely in terms of Euler characteristics. Notice, however,
that the Euler characteristics $\chi_k=\chi(\hat X \cap H_1\cap\cdots\cap H_k)$
depend not only on the variety $\hat X$ itself, but on the way in which
it is embedded into $\P^N$. This corresponds to the fact, already discussed
at length in \cite{AluMa2}, that the CSM class does not factor through the
Grothendieck ring of varieties but through a refinement, a Grothendieck
ring of immersed conical varieties.

\smallskip

The relations \eqref{CGTprimechi} and \eqref{chiandG} also give
\begin{equation}\label{chi01}
\chi(X_\Gamma) = \chi_{\hat X_\Gamma}(1) - \chi_{\hat X_\Gamma}(0) =
\sum_{k\geq 1} \chi_k(\hat X_\Gamma).
\end{equation}
In fact, we have 
$$ C'_\Gamma(T) = G'_{\mathbb{A}^n}(T) - G'_{\hat X_\Gamma}(T) = $$
$$ n(T + 1)^{n-1} - \frac{(\chi_{\hat X_\Gamma}(T + 1) + T\chi_{\hat X_\Gamma}'(T + 1))(T + 1) - T\chi_{\hat X_\Gamma}(T + 1) - \chi_{\hat{X_\Gamma}}(0)}{(T + 1)^2} $$
so that, using \eqref{CGTprimechi}, we get 
$C'_\Gamma(0) = n - (\chi_{\hat X_\Gamma}(1) - \chi_{\hat X_\Gamma}(0))=
n - \chi(X_\Gamma)$.

\medskip
\subsection{Embedded $\F_1$-structures}

The observation above leads us to consider a different, relative, 
notion of $\F_1$-structure, which takes into account not only the
variety itself but regards it as an embedded subvariety.

\smallskip

We propose the following definition of an embedded $\F_1$-structure.

\begin{defn}\label{embeddedF1}
Let $\hat X \subset \P^N$ be an embedded subvariety, defined over $\Z$,
which has the structure of an $\F_1$-variety (it admits an affine torification). 
Then $\hat X$ has an embedded
$\F_1$-structure if, in addition, the varieties $\hat X\cap H_1 \cap \cdots \cap H_k$,
where $H_1,\ldots, H_k$ are general hyperplanes, are also $\F_1$-varieties,
that is, they also admit affine torifications.
\end{defn}

\smallskip

While we focus here on torifications, the notion of embedded $\F_1$-structure 
is in general stronger than the property that $\hat X$ is an $\F_1$-variety, 
regardless of which of the currently available notions of $\F_1$-geometry one 
decides to adopt. 

\begin{rem}\label{HypRem} {\rm 
In general, one should not expect that, if $\hat X\subset \P^N$
is an $\F_1$-variety, then a hyperplane section $\hat X\cap H$ would also be.
For example, an elliptic curve (which is not a Tate motive, hence not an $\F_1$-variety)
can be realized as a hyperplane section of a $\P^2$ embedded in $\P^9$. 
However, it makes sense to look for examples of varieties satisfying this
property among those that are obtained by unions, intersections
and complements of linear subspaces, such as hyperplane arrangements. }
\end{rem}

\smallskip

The simple necessary condition on the Euler characteristic 
for the existence of an $\F_1$-structure translates immediately
into a condition expressible in terms of Chern classes for
the existence of an embedded $\F_1$-structure.

\smallskip

\begin{lem}\label{chiF1emb}
If $\hat X\subset \P^N$ has an embedded $\F_1$-structure, then
all the coefficients of the polynomial $(T+1) \, G_{\hat X}(T)$ are non-negative.
\end{lem}

\proof
The first requirement for an embedded $\F_1$ structure is that
the variety $\hat X$ is an $\F_1$-variety. This implies that 
$\hat X$ is polynomially countable and with non-negative Euler characteristic,
as we discussed in the previous sections. Thus, one obtains
the necessary condition $G_{\hat X}(0)\geq 0$.  Moreover, the requirement
that the hyperplane sections $\hat X \cap H_1 \cap \cdots \cap H_k$ are
also $\F_1$-varieties implies that the Euler characteristics
$\chi_k = \chi(\hat X \cap H_1 \cap \cdots \cap H_k)$ should also be
non-negative, and by \eqref{chiandG} this implies that all the other 
coefficients of the polynomial $(T+1) \,G_{\hat X}(T)$ should also be non-negative.
\endproof

We then obtain the following condition for the graph hypersurfaces, based
on the a positivity condition for the CSM classes.

\begin{cor}\label{posCSMXGamma}
A necessary condition for the graph hypersurface complement 
$\A^n \smallsetminus \hat X_\Gamma$ to have the structure of an
embedded $\F_1$-variety is the positivity of all the coefficients of
the polynomial $(T+1) \, C_\Gamma(T)$.
\end{cor}

\proof This follows immediately from Lemma \ref{chiF1emb} by
noticing that $$C_\Gamma(T)=G_{\A^n}(T)-G_{\hat X_\Gamma}(T)= 
G_{\A^n\smallsetminus \hat X_\Gamma}(T).$$
\endproof

\smallskip

If the coefficients of the polynomial $G_{\hat X}(T)$ are themselves non-negative,
then so are those of $(T+1) \,G_{\hat X}(T)$. Thus, effectivity of the CSM classes
implies that the necessary condition of Lemma \ref{chiF1emb} is satisfied.

\smallskip

The question of the positivity of the coefficients of the CSM class of graph
hypersurfaces was first raised in \cite{AluMa1}, where it is conjectured that
this condition may always be satisfied. However, Aluffi has recently shown 
that there are examples of graph hypersurfaces, for sufficiently
large graphs, where the positivity condition may fail, \cite{AluPC}. Thus, the necessary
condition for  the existence of embedded $\F_1$-structure may not be always satisfied
for graph hypersurfaces.

\smallskip

It is important to stress the fact that the condition of Lemma \ref{chiF1emb}
is certainly not sufficient, even when combined with the other necessary
condition that the varieties $\hat X \cap H_1 \cap\cdots \cap H_k$ are 
polynomially countable (Tate motives). For example, it is known from the
results of \cite{AluMi}, \cite{Jones}, \cite{Stry}, that the CSM classes of many 
Schubert varieties satisfy the positivity condition of Lemma \ref{chiF1emb},
but Schubert varieties are not necessarily $\F_1$-varieties: for example,
in the torified-schemes formulation of $\F_1$-geometry developed in \cite{LL},
Grassmannians and Schubert varieties are torified  varieties but not necessarily 
endowed with an affine torification, see the discussion in Example 1.15 of \cite{LL}. 
We will return to discuss the role of torifications in \S \ref{TorSec}.

\medskip
\subsection{A $q$-deformed CSM class}

Let $\hat X \subset \P^N$, defined over $\Z$, be a polynomially countable subvariety with
the property that its hyperplane sections are also polynomially countable. 
This condition allows for the definition of a $q$-deformed CSM-class,
whose limit as $q\to 1$ gives back the invariant $G_{\hat X}(T)$.

\begin{defn}\label{defCSMq}
Let $\hat X \subset \P^N$, defined over $\Z$, be such that $\hat X \cap H_1 \cap \cdots \cap H_k$ is
polynomially countable, for general hyperplanes $H_j$. Then the $q$-deformed
CSM class is given by
\begin{equation}\label{CSMq}
G_{\hat X}(q,T) := \frac{T N_{\hat X}(q,T) + N_{\hat X}(q)}{T+1},
\end{equation}
where the polynomial $N_{\hat X}(q,T)$ is defined as
\begin{equation}\label{NXqT}
N_{\hat X}(q,T) =\sum_{k\geq 0} N_{\hat X\cap H_1 \cap\cdots \cap H_k}(q) T^k,
\end{equation}
with $N_{\hat X\cap H_1 \cap\cdots \cap H_k}(q)$ the polynomial interpolating
the number of points over $\F_q$ of the mod $p$ reduction of $\hat X\cap H_1 \cap
\cdots \cap H_k$, and $N_{\hat X}(q)=N_{\hat X}(q,0)=\# \hat X_p(\F_q)$, with 
$\hat X_p$ the mod $p$ reduction of $\hat X$.
\end{defn}

While the counting of points over finite fields is defined for $q=p^r$ a prime power,
the polynomially countable hypothesis implies that all the 
$N_{\hat X\cap H_1 \cap\cdots \cap H_k}(q)$ are polynomials,
defined for any real or complex value of $q$. Thus, it makes sense to
take a limit when $q\to 1$.
By construction, this limit of $G_{\hat X}(q,T)$ gives back
the usual CSM class, in the form of the polynomial $G_{\hat X}(T)$.

\medskip
\section{Regular torifications and Chern classes}\label{TorSec}

In \cite{Alu2} an explicit way is obtained to compute a generalization 
of the CSM class $\csm(X)$ for any variety $X$, which reduces to the 
usual CSM class when $X$ is complete. 
This is called the pro-Chern class, and is also denoted by $\csm$. It can be
constructed in terms of a notion of good closure and a pro-Chow functor. 
Namely, let $U$ be a smooth variety with closure $\overline{U}$ such that $\overline{U}$ is a \textit{good} closure, that is smooth with $\overline{U} \setminus U$ a normal crossings divisor. Let $\{U\} = c(\Omega_{\overline{U}}^1(\overline{U}\setminus U)^\vee) \cap [\overline{U}]$ where $c$ denotes the usual Chern class of a bundle and $\Omega_{\overline{U}}^1(\overline{U}\setminus U)$ are the differentials with logarithmic poles on $\overline{U} \setminus U$. Then $\csm(U) = \{U\}$. Let $X$ be any variety and suppose $\mbf{1}_X = \sum m_U \mbf{1}_U$ where $U$ are smooth, locally closed subvarieties and $m_U$ are integers. Then the class $\csm(X)$ is
expressed as a finite sum 
$\csm(X) = \sum m_U {i_U}_* \{U\}$, where $i_U: U \to X$ is the inclusion of $U$ into $X$ and more generally $\csm(\alpha) = \sum m_U {i_U}_* \{U\}$, whenever $\alpha = \sum m_U \mbf{1}_U$ is a constructible function on $X$. Theorem 3.3 of \cite{Alu2} shows that this is independent of compatible good closures. That is, if $\overline{U}$ is the closure of $U$ in $X$ with inclusion $i$ and $\overline{U}'$ is a good closure of $U$ with inclusion $j$, such that there exists a morphism $\pi : \overline{U} \to \overline{U}'$ with $j = \pi \circ i$, then both these closures give the same class $\{U\}$ up to passing to the pro-Chow  homology of $X$. 

\medskip
\subsection{Chern class for toric varieties and positivity}

An interesting application, given in \cite{Alu2}, of this construction of CSM classes
is a new simple derivation of the Ehlers formula for the Chern classes of toric varieties,
which in particular shows that positivity holds, since the Chern class reduces to a
sum of fundamental classes of closures of torus orbits.

More precisely, it is proved in Theorem 4.2 of \cite{Alu2} that the pro-CSM class
of a toric variety $X$ can be written as
\begin{equation}\label{csmXtoric}
\csm(X) = \sum_{O \in X/T} c(\Omega^1_{\bar O}(\log D)^\vee) \cap [\bar O] 
= \sum_{O \in X/T} [\bar O],
\end{equation}
where $O$ ranges over the orbits of the action of the torus $T$ on $X$, with
$[\bar O]$ the fundamental class of the orbit closure in the pro-Chow homology of $X$,
and where $D=\bar O\smallsetminus O$. The second equality depends on the well
known fact (\cite{Ful}, p.87) that $\Omega^1_{\bar O}(\log D)$ is trivial in this case. 

\medskip
\subsection{Regular torifications}

In general, one does not expect that a similar result would hold for torified
varieties. In fact, in a torification in general there is no torus action, unlike in
the case of the toric varieties, and the natural torus action one has on each of the
tori $T_i$ of the torification in $X$ does not in general extend to the closure
$\bar T_i$ in $X$.

However, there is a better behaved class of torifications, which were
considered in \S 6.2 of \cite{LL} as a possible way to eliminate the
ambiguity, due to the choice of torification, that can lead to different $\F_1$-structures
on the same $\Z$-manifold. These are the {\em regular} torifications.

\begin{defn}\label{regtordef} {\rm (L\'opez Pe\~{n}a and Lorscheid \cite{LL})} 
A torification $e_X : T=\amalg_{i\in I} T_i \to X$ is regular if, for all tori $T_i\subset T$, 
there exists a set $J_i\subseteq I$ such that the closure $\bar T_i$ in $X$ is a union
of tori in the same torification, $\bar T_i =\amalg_{j\in J_i} T_j$.
\end{defn}

It is not known, at present, whether all varieties admitting a torification
(or an affine torification) necessarily admit a regular one, though this
is the case in all the explicit examples studied so far.

\smallskip

Thus, it is interesting to ask whether the argument of Theorem 4.2 of \cite{Alu2}
for the computation of CSM classes would extend to the case of regular
torifications. Notice that, because all the flag varieties admit a regular torification,
this would then automatically prove the conjectured positivity of the CSM classes
for such varieties. The question is therefore whether, given a regular
torification $e_X : T \to X$, the pro-Chern class
$$ \csm(X) = \sum_{T_i \in T} c(\Omega^1_{\bar T_i}(\log D_i)^\vee) \cap [\bar T_i] $$
can be written as
$$  \csm(X) = \sum_{T_i \in T}  [\bar T_i]. $$
This would be the case if we knew that $\Omega^1_{\bar T_i}(\log D_i)$ is trivial. 

\smallskip

However, this will in general not be the case. 

\begin{lem}\label{torDlem}
Let $e: T=\bG_m^d \hookrightarrow X$ be an embedding of a torus in a smooth
projective variety $X$ such that $D=\bar T \smallsetminus T$ is a 
simple normal crossings divisor, but  $\bar T$ is not an equivariant
compactification of a semi-abelian variety. Then 
$\Omega^1_{\bar T}(\log D)$ is non-trivial. 
\end{lem}

\proof
There is a characterization (\cite{Wink1}, \cite{Wink2}) of varieties with
trivial logarithmic tangent bundle. In particular, by Corollary 1 of \cite{Wink1},
if $X$ is a non-singular complete algebraic complex variety and $D\subset X$
a simple normal crossings divisor, then the condition that the logarithmic tangent
bundle $\Omega^1_X(\log D)^\vee$ is trivial is equivalent to the existence of a semi-abelian
variety $A$ that acts on $X$ with $X\smallsetminus D$ an open orbit.
\endproof

Therefore, one can construct examples where the argument of Theorem 4.2 of \cite{Alu2}
does not directly extend to regular torifications, by producing an example of a regular
torification $e_X: T \to X$, where some of the tori $T_i$ satisfy the condition of
Lemma \ref{torDlem}.

\smallskip

However, there are cases where one may still be able to use a {\em regular} 
torification to obtain useful information on the CSM classes. For instance, 
consider the case of flag varieties. It is known (see \S\S 1.3.5 and 6.2 of \cite{LL}) 
that these admit a unique isomorphism class of regular torifications, coming from
Schubert cell decomposition (these are in general not affine torifications, though). 
In \cite{AluMi} an explicit construction of the Chern
classes of Schubert varieties is given in terms of the decomposition in Schubert
cells, and a computation of the CSM classes of Schubert cells in terms of certain
uniquely determined integer coefficients (which are conjectured to be always
non-negative). In view of the considerations above, one can ask the following  
question.

\begin{ques}\label{Schubertques}
Is there a description of the coefficients $\gamma_{\underline{\alpha},\underline{\beta}}$
in the computation of the CSM classes of Schubert varieties in \cite{AluMi} in terms of
counting of tori in a (unique) regular torification? Can such an approach 
be used to prove positivity?
\end{ques}

\bigskip

\subsection*{Acknowledgment} The first author was supported for this project by
a Caltech Summer Undergraduate Research Fellowship. The second author
is supported by NSF grants DMS-0901221, DMS-1007207, DMS-1201512, 
PHY-1205440. The second author thanks Paolo Aluffi for many useful discussions
and for a careful reading of the manuscript, Javier L\'opez-Pe\~{n}a for reading
an earlier draft of the paper and offering comments and suggestions, and 
Spencer Bloch for useful conversations about \cite{Bloch}.

\end{document}